\documentclass[twocolumn]{autart}

\input{latex_common/preamble}

\newcommand{\union}{\cup}
\newcommand{\intersection}{\cap}

\newcommand{\set}[1]{\mathcal #1}

\newcommand{\reals}{\mathbb R}
\newcommand{\real}{\mathbb R}

\newcommand{\interior}[1]{\mathrm{int}(#1)}

\newcommand{\boundary}[1]{\mathrm{bdry}\left(#1\right)}

\newcommand{\algvar}[1]{\text{\IfSubStr{#1}{_}{%
    \StrSubstitute{#1}{_}{\textunderscore}}{#1}}}

\newcommand{\grad}{\nabla}

\newcommand{\subdiff}[1][none]{%
  \ifthenelse{\equal{#1}{none}}{%
    \partial%
  }{%
    \partial_{#1}%
  }%
}
\newcommand{\normalcone}[2]{\set{N}_{#1}(#2)}
\DeclareMathOperator{\conehull}{cone}

\DeclareMathOperator{\range}{range}

\newcommand{\diag}[2][short]{\mathrm{diag}\ifthenelse{\equal{#1}{short}}{(#2)}{\left(#2\right)}}

\newcommand{\norm}[2]{\|#1\|_{#2}}
\newcommand{\vertiii}[1]{{\left\vert\kern-0.25ex\left\vert\kern-0.25ex\left\vert #1 
        \right\vert\kern-0.25ex\right\vert\kern-0.25ex\right\vert}}

\newcommand{\Proj}[2]{\mathcal P_{#1}(#2)}

\newcommand{\PP}{$\mathcal{P}$\xspace}

\newcommand{\NPhard}{$\mathcal{NP}$-hard\xspace}

\newcommand{\transp}{{\scriptscriptstyle\mathsf{T}}}
\newcommand{\T}{^\transp}

\newcommand{\fun}[2][1]{%
  #2(%
  \foreach \index in {1, ..., #1} {%
    \ifthenelse{\equal{\index}{#1}}{%
      \cdot%
    }{%
      \cdot,%
    }%
  })}
\newcommand{\definedas}[1][tri]{%
  \ifthenelse{\equal{#1}{tri}}{\triangleq}{\coloneqq}}

\newcommand{\inlinesum}{\textstyle\sum}

\DeclareMathOperator*{\exptx}{exp}
\renewcommand{\exp}[2][exponent]{\ifthenelse{\equal{#1}{exponent}}{e^{#2}}{\exptx\left(#2\right)}}
\DeclareMathOperator*{\expsf}{\mathsf{exp}}
\newcommand{\expm}[2][exponent]{\ifthenelse{\equal{#1}{exponent}}{\mathsf{e}^{#2}}{\expsf\left(#2\right)}}

\makeatletter
\DeclareFontFamily{U}{tipa}{}
\DeclareFontShape{U}{tipa}{m}{n}{<->tipa10}{}
\newcommand{\arc@char}{{\usefont{U}{tipa}{m}{n}\symbol{62}}}%
\renewcommand{\arc}[1]{\mathpalette\arc@arc{#1}}
\newcommand{\arc@arc}[2]{%
  \sbox0{$\m@th#1#2$}%
  \vbox{
    \hbox{\resizebox{\wd0}{\height}{\arc@char}}
    \nointerlineskip
    \box0
  }%
}
\makeatother

\newcommand{\alev}[1]{\text{a.e. }#1}

\newcommand{\ev}[3][c]{%
  \ifthenelse{\equal{#1}{c}}{%
    \ifthenelse{\equal{#2}{}}{\forall[0,#3]}{\forall[#2,#3]}%
  }{%
    \ifthenelse{\equal{#1}{o}}{%
      \ifthenelse{\equal{#2}{}}{\forall(0,#3)}{\forall(#2,#3)}%
    }{%
      \ifthenelse{\equal{#1}{oc}}{%
        \ifthenelse{\equal{#2}{}}{\forall(0,#3]}{\forall(#2,#3]}%
      }{%
        \ifthenelse{\equal{#2}{}}{\forall[0,#3)}{\forall[#2,#3)}%
      }%
    }%
  }%
}

\DeclareMathOperator*{\argmax}{argmax}
\DeclareMathOperator*{\argmin}{argmin}

\def\opticmd{\min}

\newcommand{\optimal}[1]{#1^*}

\newcounter{l}
\newcounter{j}
\newcounter{k}

\makeatletter
\newenvironment{taggedsubequations}[1]
{%
  \addtocounter{equation}{-1}%
  \begin{subequations}%
    \def\@currentlabel{#1}%
    \renewcommand{\theequation}{#1.\alph{equation}}%
  }
  {\end{subequations}}
\makeatother 

\newcommand{\lloptimization}[7][center]{
  \setsepchar{\#}%
  \readlist\mylist{#6}
  \ifthenelse{\equal{#7}{}}{\begin{subequations}}{\begin{taggedsubequations}{#7}}
    \ifthenelse{\equal{#2}{}}{}{\label{eq:#2}}
    \ifthenelse{\equal{#1}{left}}{
      \begin{flalign}
        \ifthenelse{\equal{#2}{}}{}{\label{eq:#2_a}}
        \ifthenelse{\equal{#3}{}}{}{#3 = }
        &\opticmd_{\ifthenelse{\equal{#4}{}}{}{#4}}~ #5~\mathrm{s.t.}\hspace{-1mm} &&\\
        \forloop{k}{0}{\arabic{k} < \listlen\mylist[]}{
          \setcounter{j}{\value{k}+1}
          \setcounter{l}{\value{k}+2}
          \ifthenelse{\equal{#2}{}}{}{\label{eq:#2_\alph{l}}}
          \ifthenelse{\equal{\arabic{j}}{\listlen\mylist[]}}{%
            &\mylist[\arabic{j}]%
          }{%
            &\mylist[\arabic{j}] &&\\%
          }%
        }
      \end{flalign}
    }{\ifthenelse{\equal{#1}{center}}{
      \begin{align}
        \ifthenelse{\equal{#2}{}}{}{\label{eq:#2_a}}
        \ifthenelse{\equal{#3}{}}{}{#3 = }
        &\opticmd_{\ifthenelse{\equal{#4}{}}{}{#4}}~ #5~\mathrm{s.t.}\hspace{-1mm} &&\\
        \forloop{k}{0}{\arabic{k} < \listlen\mylist[]}{
          \setcounter{j}{\value{k}+1}
          \setcounter{l}{\value{k}+2}
          \ifthenelse{\equal{#2}{}}{}{\label{eq:#2_\alph{l}}}
          \ifthenelse{\equal{\arabic{j}}{\listlen\mylist[]}}{%
            &\mylist[\arabic{j}]%
          }{%
            &\mylist[\arabic{j}]\\%
          }%
        }
      \end{align}}{\ifthenelse{\equal{#1}{left*}}{
      \begin{flalign*}
        \ifthenelse{\equal{#3}{}}{}{#3 = }
        &\opticmd_{\ifthenelse{\equal{#4}{}}{}{#4}}~ #5~\mathrm{s.t.}\hspace{-1mm} &&\\
        \forloop{k}{0}{\arabic{k} < \listlen\mylist[]}{
          \setcounter{j}{\value{k}+1}
          \setcounter{l}{\value{k}+2}
          \ifthenelse{\equal{\arabic{j}}{\listlen\mylist[]}}{%
            &\mylist[\arabic{j}]%
          }{%
            &\mylist[\arabic{j}] &&\\%
          }%
        }
      \end{flalign*}
    }{\begin{align*}
        \ifthenelse{\equal{#3}{}}{}{#3 = }
        &\opticmd_{\ifthenelse{\equal{#4}{}}{}{#4}}~ #5~\mathrm{s.t.}\hspace{-1mm} &&\\
        \forloop{k}{0}{\arabic{k} < \listlen\mylist[]}{
          \setcounter{j}{\value{k}+1}
          \setcounter{l}{\value{k}+2}
          \ifthenelse{\equal{\arabic{j}}{\listlen\mylist[]}}{%
            &\mylist[\arabic{j}]%
          }{%
            &\mylist[\arabic{j}]\\%
          }%
        }
      \end{align*}}}
    }
  \ifthenelse{\equal{#7}{}}{\end{subequations}}{\end{taggedsubequations}}
}

\newcommand{\llpoptimization}[8][left]{
  \begin{problem}[#8]\ifthenelse{\equal{#2}{NoNe}}{}{\textbf{#2.}}%
    \ifthenelse{\equal{#3}{}}{}{\label{problem:#3}}%
    \lloptimization[#1]{#3}{#4}{#5}{#6}{#7}{#8}%
  \end{problem}%
}

\NewEnviron{optimization}[6]{\lloptimization[#1]{#2}{#3}{#4}{#5}{\BODY}{#6}}
\NewEnviron{foptimization}[7][NoNe]{%
  \begin{mdframed}[default,linewidth=1pt]%
    \llpoptimization[#2]{#1}{#3}{#4}{#5}{#6}{\BODY}{#7}%
  \end{mdframed}%
}
\NewEnviron{poptimization}[7][NoNe]{%
  \llpoptimization[#2]{#1}{#3}{#4}{#5}{#6}{\BODY}{#7}%
}

\definecolor{darkolivegreen}{rgb}{0.33, 0.6, 0.18}
\definecolor{green}{rgb}{0, 0.5, 0}
\definecolor{orange}{rgb}{1, 0.4, 0}
\definecolor{darkgray}{rgb}{0.2, 0.2, 0.2}

\newcommand{\postmeeting}[2][show]{%
  \ifthenelse{\equal{#1}{show}}{%
    {\color{orange} #2}%
  }{}%
}
\newcommand{\idea}[2][show]{%
  \ifthenelse{\equal{#1}{show}}{%
    {\color{orange} #2}%
  }{}%
}
\newcommand{\todo}[2][show]{%
  \ifthenelse{\equal{#1}{show}}{%
    {\color{blue}%
      \ifthenelse{\equal{#2}{}}{TODO}{(TODO: #2)}}%
  }{}%
}
\newcommand{\question}[2][show]{%
  \ifthenelse{\equal{#1}{show}}{%
    {\color{darkolivegreen} (Q: #2)}%
  }{}%
}
\newcommand{\fixme}[2][show]{%
  \ifthenelse{\equal{#1}{show}}{%
    {\color{red} (FIXME: #2)}%
  }{}%
}
\newcommand{\mycomment}[2][show]{%
  \ifthenelse{\equal{#1}{show}}{%
    {\color{red} (C: #2)}%
  }{}%
}

\makeatletter
\@ifclassloaded{beamer}{}{}
\makeatother

\renewcommand{\algref}[3]{\ifthenelse{\equal{#2}{}}{Algorithm~\ref{alg:#1}\xspace}%
  {\ifthenelse{\equal{#3}{}}{line~\ref{alg:#1:line:#2} %
      of Algorithm~\ref{alg:#1}\xspace}{lines~%
      \ref{alg:#1:line:#2}-\ref{alg:#1:line:#3} %
      of Algorithm~\ref{alg:#1}\xspace}}}

\allowdisplaybreaks

\tcbset{highlight math style={enhanced,
colframe=red!60!black,colback=yellow!50!white,arc=4pt,boxrule=1pt,
drop fuzzy shadow}}

\newcolumntype{C}{>{\centering\arraybackslash}X}
\newcolumntype{L}{>{\raggedright\arraybackslash}X}
\newcolumntype{R}{>{\raggedleft\arraybackslash}X}

\makeatletter
\@ifclassloaded{beamer}{}{
\@ifclassloaded{autart}{}{

\theoremstyle{definition}

\theoremstyle{empty}

\declaretheoremstyle[
  headfont=\color{black}\normalfont\bfseries,
  bodyfont=\color{darkgray}\normalfont,
]{colored}

}
}
\makeatother

\newcounter{theorem}
\newtheorem{innertheorem}{Theorem}
\makeatletter
\newenvironment{theorem}[1]
{\refstepcounter{theorem}%
\protected@edef\@currentlabelname{\value{theorem}}%
\renewcommand\theinnertheorem{{\thetheorem#1}}
\innertheorem}
{\endinnertheorem}
\makeatother

\newtheorem{lemma}{Lemma}

\theoremstyle{definition}

\newtheorem{definition}{Definition}

\newcounter{condition}
\newtheorem{innercondition}{Condition}
\makeatletter
\newenvironment{condition}[2][]
{\refstepcounter{condition}%
\protected@edef\@currentlabelname{\value{condition}}%
\ifthenelse{\equal{#2}{}}{\renewcommand\theinnercondition{{\thecondition#1}}}
{\renewcommand\theinnercondition{{\thecondition#1 \textnormal{(#2)}}}}
\innercondition}
{\endinnercondition}
\makeatother

\newcounter{problem}

\newtheorem{innerproblem}{Problem}
\makeatletter
\newenvironment{problem}[1][]
{\ifthenelse{\equal{#1}{}}
{\refstepcounter{problem}%
\protected@edef\@currentlabelname{\value{problem}}%
\renewcommand\theinnerproblem{\theproblem}\innerproblem}
{\renewcommand\theinnerproblem{#1}\innerproblem}}
{\endinnerproblem}
\makeatother

\newtheorem{assumption}{Assumption}

\newmdenv[
outerlinewidth = 1,%
linewidth = 0pt,%
roundcorner = 2pt,%
leftmargin = 20,%
rightmargin = 0,%
backgroundcolor = lightgray!10,%
outerlinecolor = gray!50,%
innertopmargin = \topskip,%
splittopskip = \topskip,%
frametitle = Summary,%
frametitlebelowskip = 0pt,%
]{summary_box}

\newmdenv[
outerlinewidth = 1,%
linewidth = 0pt,%
roundcorner = 2pt,%
leftmargin = 60,%
rightmargin = 60,%
backgroundcolor = yellow!40,%
outerlinecolor = black!50,%
innertopmargin = \topskip,%
splittopskip = \topskip,%
]{highlight_box}

\newif\ifnomentry

\renewcommand{\nomgroup}[1]{%
\nomentryfalse
\ifthenelse{\equal{#1}{V}}{\item[\textbf{Variables}]}{%
\ifthenelse{\equal{#1}{A}}{\item[\textbf{Abbreviations}]}{%
\ifthenelse{\equal{#1}{E}}{\item[\textbf{Equations}]}{}}}
\nomentrytrue
}
\newcommand{\defvar}[3][show]{\nomenclature[V]{#2}{#3}\ifthenelse{\equal{#1}{show}}{#2}{}}

\algrenewcommand\algorithmicindent{1em}
\algrenewcommand\alglinenumber[1]{{\sf\footnotesize#1}}
\algtext*{EndWhile}
\algtext*{EndIf}
\algtext*{EndFor}

\algrenewcommand\algorithmicrequire{\textbf{Precondition:}}
\algrenewcommand\algorithmicensure{\textbf{Postcondition:}}
\algnewcommand\algorithmicinput{\textbf{Input:}}
\algnewcommand\Input{\item[\algorithmicinput]}
\algnewcommand\algorithmicfunctions{\textbf{Lambda Functions:}}
\algnewcommand\Functions{\item[\algorithmicfunctions]}
\algnewcommand\Stop[1][]{\State \textbf{STOP} {\color{green} \ul{#1}}}
\algnewcommand\Break{\State \textbf{break}}
\algnewcommand\Continue{\State \textbf{continue}\xspace}
\algnewcommand\Error[1][]{\State \textbf{error} {\color{red} \uwave{#1}}}
\algrenewcommand\algorithmiccomment[1]{\hfill {\color{gray} \(\triangleright\) #1}}
\makeatletter
\algnewcommand\Section[1]{\Statex \hskip\ALG@thistlm \textbf{\color{gray} \(\triangleright\) #1}}
\makeatother
\algdef{SE}[PARAMETERS]{Parameters}{EndParameters}
   {\algorithmicparameters}
   {\algorithmicend\ \algorithmicparameters}
\algnewcommand{\algorithmicparameters}{\textbf{parameters}}
\newcommand{\algorithmiclambda}[3][]{\ifthenelse{\equal{#1}{}}{\texttt{#2}(#3)}{$\texttt{#2}\gets\textbf{lambda}~#3:~#1$}}

\newcommand{\tzero}{0}

\usepackage[T1]{fontenc}
\usepackage[utf8]{inputenc}

\graphicspath{{./figures/}}

\begin{document}

\begin{frontmatter}

  \title{Lossless convexification of non-convex optimal control problems with
    disjoint semi-continuous inputs\thanksref{footnoteinfo}}

  \thanks[footnoteinfo]{This paper was not presented at any IFAC
    meeting. Corresponding author D.~Malyuta.}

  \author[Student]{Danylo Malyuta}\ead{danylo@uw.edu},
  \author[Student]{Michael Szmuk}\ead{michael.szmuk@gmail.com},
  \author[Professor]{Beh\c{c}et A\c{c}{\i}kme\c{s}e}\ead{behcet@uw.edu}

  \address[Student]{Doctoral Student, Dept. of Aeronautics \& Astronautics,
    University of Washington, Seattle, WA 98195, USA}
  \address[Professor]{Professor, Dept. of Aeronautics \& Astronautics,
    University of Washington, Seattle, WA 98195, USA}


  \begin{abstract}
    This paper presents a convex optimization-based method for finding the
    globally optimal solutions of a class of mixed-integer non-convex optimal
    control problems. We consider problems that are non-convex in the input
    norm, which is a semi-continuous variable that can be zero or lower- and
    upper-bounded. Using lossless convexification, the non-convex problem is
    relaxed to a convex problem whose optimal solution is proved to be
    optimal almost everywhere for the original problem. The relaxed problem can be
    solved using second-order cone programming, which is a subclass of convex
    optimization for which there exist numerically reliable solvers with
    convergence guarantees and polynomial time complexity. This is the first
    lossless convexification result for mixed-integer optimization problems. An
    example of spacecraft docking with a rotating space station corroborates the
    effectiveness of the approach and features a computation time almost three
    orders of magnitude shorter than a mixed-integer programming formulation.
  \end{abstract}  
\end{frontmatter}

\section{Introduction}

We present a convex programming solution to a class of optimal control problems
with semi-continuous control input norms. Semi-continuous variables are a
particular type of binary non-convexity.

\begin{definition}
  \label{definition:semicontinuous}
  Variable $x\in\reals$ is \textnormal{semi-continuous} if
  $x\in \{0\}\union[a,b]$ with $0<a\le b$ \cite{MosekCookbook2019}.
\end{definition}

The constraint $az\le x\le bz$ with $z\in\{0,1\}$ models semi-continuity. We
consider systems that have multiple inputs which may not all be simultaneously
active, which point in dissimilar directions in the input space, and whose norms
are semi-continuous. Although mixed-integer convex programming (MICP) is
applicable, it is an \NPhard optimization class \cite{Bemporad1999,Cormen2009}
and solving a practical path planning problem such as rocket landing or
spacecraft rendezvous can take hours. This paper proposes an algorithm based on
lossless convexification that solves these problems to global optimality in
seconds.

Non-convex lower-bound constraints on the input norm have been handled in past
research using convex optimization via lossless convexification. In this method
the original problem is relaxed to a convex one via a slack variable, enabling
the use of second-order cone programming (SOCP) to solve the original problem to
global optimality in polynomial time. The method was introduced in
\cite{Acikmese2007} for minimum-fuel rocket landing and was later expanded to
fairly general non-convex input sets \cite{Acikmese2011}. Extensions of the
method were introduced in \cite{Blackmore2010,
  Carson2011,Acikmese2013} to handle minimum-error rocket landing and non-convex
pointing constraints. The method was used in \cite{Pascucci2017} for satellite
docking trajectory generation. More recently, lossless convexification was
rigorously shown to handle affine and quadratic state constraints
\cite{Harris2013a,Harris2013b}, culminating in \cite{Harris2014} which has
to-date been the most general formulation. However, a recurring assumption is
that there is a single input which cannot be turned off.

Our interest is in problems with multiple such inputs, which are allowed to turn
off and which may not all be simultaneously active. Such problems can be solved
directly by wrapping existing lossless convexification results in a
mixed-integer program. Binary variables would decide which inputs are
active. Similar ideas were explored in \cite{Blackmore2012,Zhang2017}. However,
the \NPhard nature of MICP makes the approach computationally expensive and
without a real-time guarantee.

An alternative solution is through successive convexification, where
non-linearities are linearized and a sequence of convex programs is solved until
convergence \cite{Mao2016,Mao2017,Dueri2017,Mao2018,Bonalli2019}. Lossless
convexification can similarly be embedded to handle input non-convexity
\cite{Szmuk2016}. Although the method was originally devised for optimal control
problems with continuous variables, an effective way has recently been found to
embed binary decisions in a continuous formulation
\cite{Szmuk2019a,Szmuk2019b,Reynolds2019,Malyuta2019}. Although successive
convexification is faster than mixed-integer programming, it is a local
optimizer and is not guaranteed to converge to a feasible solution.

Our main contribution is to extend lossless convexification to directly handle a
class of mixed-integer non-convex optimal control problems with multiple inputs
and semi-continuous input norms in the sense of
Definition~\ref{definition:semicontinuous}. Unlike mixed-integer programming,
lossless convexification solves the problem in polynomial time. Unlike
successive convexification, it finds the global optimum and is guaranteed to
converge. 
The approach is amenable to real-time
onboard optimization for autonomous systems or for rapid design trade studies.

The paper is organized as follows. Section~\ref{sec:problem_statement} defines
the class of optimal control problems that our method
handles. Section~\ref{sec:lossless_convexification} then proposes our solution
method. Section~\ref{sec:ct_lcvx_proof} proves that our method finds the
globally optimal solution based on the necessary conditions of optimality
presented in Section~\ref{sec:pmp}. Section~\ref{sec:examples} presents an
example which corroborates the method's effectiveness for practical control
problems. Section~\ref{sec:future_work} outlines future work and
Section~\ref{sec:conclusion} summarizes the result.

\textit{Notation}: sets are calligraphic, e.g. $\mathcal S$. Operator $\circ$
denotes the element-wise product. Given a function
$f:\reals^n\times\reals^m\to\reals^p$, we use the shorthand
$f[t]\equiv f(x(t),y(t))$. In text, functions are referred to by their letter
(e.g. $f$) and conflicts with another variable are to be understood from
context. The gradient of $f$ with respect an argument $x$ is denoted
$\grad_x f\in\reals^{p\times n}$. Similarly, if $f$ is nonsmooth then its
subdifferential with respect to $x$ is
$\subdiff[x] f\subseteq\reals^{1\times n}$. The normal cone at $x$ to
$\mathcal S\subseteq\reals^n$ is denoted
$\normalcone{\mathcal S}{x}\subseteq\reals^n$. Given vectors $x_i\in\reals^n$,
$\conehull\{x_1,x_2,\dots\}\subseteq\reals^n$ denotes their conical hull. When
we refer to an \textit{interval}, we mean some time interval $[t_1,t_2]$ of
non-zero duration, i.e. $t_1<t_2$. We
call the Eucledian projection of $y\in\reals^n$ onto $\set S\subseteq\reals^n$
the magnitude of the 2-norm projection of $y$:
\begin{eqnarray}
  \label{eq:eucledian_projection}
  \Proj{\set S}{y}\definedas\big\|\textstyle\argmin_{z\in\set S}\norm{y-z}{2}\big\|_2.
\end{eqnarray}

\section{Problem Statement}
\label{sec:problem_statement}

We consider mixed-integer non-convex optimal control problems that extend the
problem class defined in \cite{Acikmese2011}:
\begin{poptimization}{left}{ocp}{}{u_i,\gamma_i,t_f}{\hspace{-1mm}m(t_f,x(t_f))}
  {$\mathcal O$}
  \dot x(t) = Ax(t)+B\inlinesum_{i=1}^M{u_i(t)}+w,~x(0) = x_0,\hspace{-1mm} \#
  \gamma_i(t)\rho_1 \le \|u_i(t)\|_2\le \gamma_i(t)\rho_2\quad i=1,\dots,M, \#
  \gamma_i(t)\in \{0,1\}\quad i=1,\dots,M, \#
  \inlinesum_{i=1}^M\gamma_i(t) \le K, \#
  C_iu_i(t)\le 0\quad i=1,\dots,M, \#
  b(t_f,x(t_f)) = 0,
\end{poptimization}
where $x(t)\in\reals^n$ is the state, $u_i(t)\in\reals^m$ is the $i$-th input,
and $w\in\reals^n$ is a known external input. Convex functions
$m:\reals\times\reals^n\to\reals$ and $b:\reals\times\reals^n\to\reals^{n_b}$
define the terminal cost and the terminal manifold respectively. The input
directions are constrained to polytopic cones called \textit{input pointing
  sets}:
\begin{eqnarray}
  \label{eq:input_pointing_set}
  \mathcal U_i\definedas \{u\in\reals^m: C_iu\le 0\},
\end{eqnarray}
where $C_i\in\reals^{p_i\times m}$ is a matrix with $C_{i,j}$ the $j$-th row,
such that $C_{i,j}\T$ defines the outward-facing normal of the $j$-th
facet. The constraint \eqref{eq:ocp_f} allows at most $K\le M$ inputs to be
turned on simultaneously.

\begin{assumption}
  \label{ass:pointing_set_non_overlapping}
  The pointing set interiors do not overlap,
  i.e. $\interior{\mathcal U_i}\cap\interior{\mathcal U_j}=\emptyset$ for all
  $i\ne j$.
\end{assumption}

\begin{assumption}
  \label{ass:full_range_cone}

  Matrices $C_i$ in \eqref{eq:ocp_f} are full row rank and the terminal cost is
  non-trivial, i.e. $\grad m\ne 0$.
\end{assumption}

\begin{assumption}
  \label{ass:thrust_bounds}
  The control norm bounds in \eqref{eq:ocp_c} are distinct,
  i.e. $\rho_1<\rho_2$.
\end{assumption}

\section{Lossless Convexification}
\label{sec:lossless_convexification}

\begin{figure*}
  \centering
  \begin{subfigure}[b]{0.32\textwidth}
    \centering
    \includegraphics[width=0.87\columnwidth]{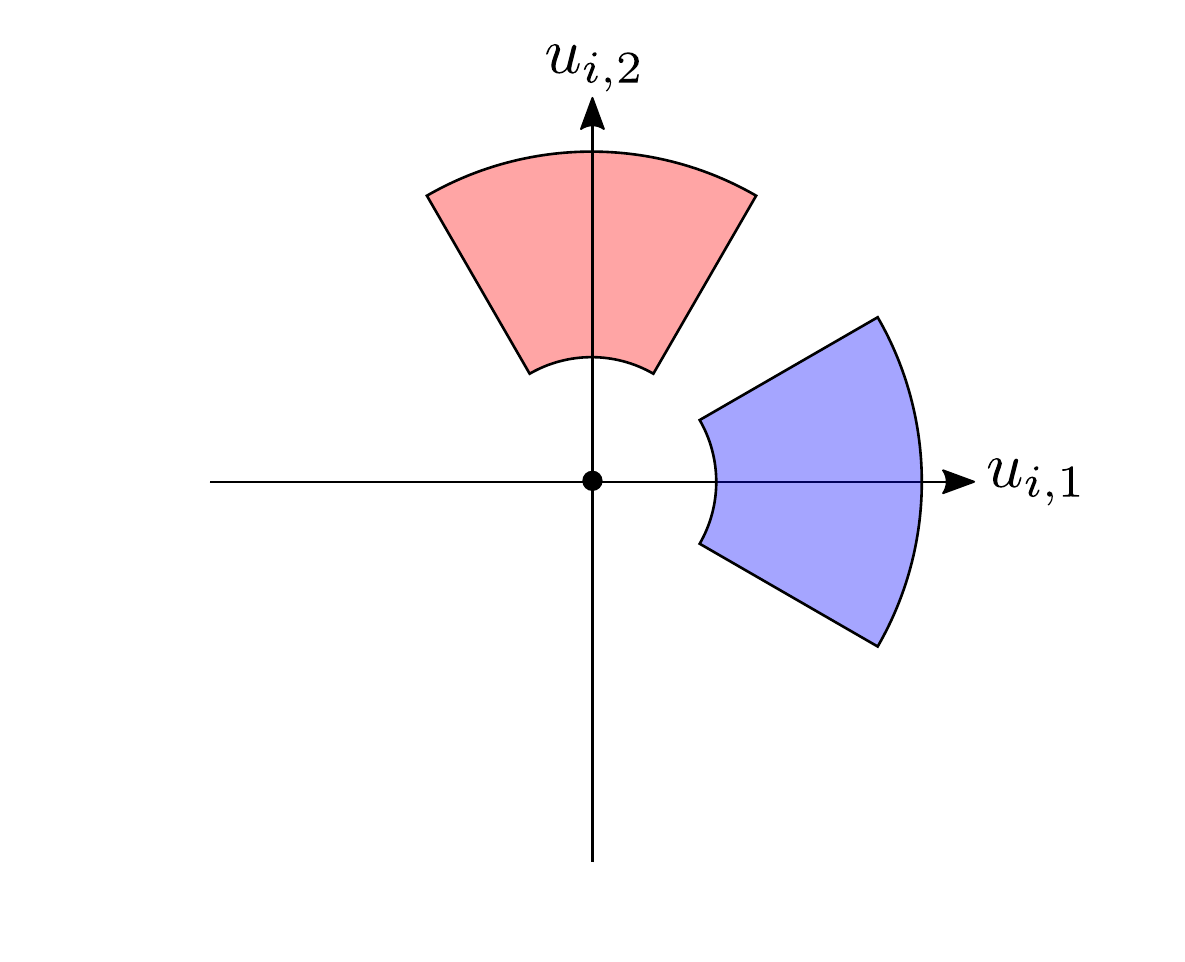}
    \caption{Original non-convex disjoint input sets defined by
      \eqref{eq:ocp_c}-\eqref{eq:ocp_f}.}
    \label{fig:input_set_relaxation_1}
  \end{subfigure}%
  \hfill%
  \begin{subfigure}[b]{0.32\textwidth}
    \centering
    \includegraphics[width=0.9\columnwidth]{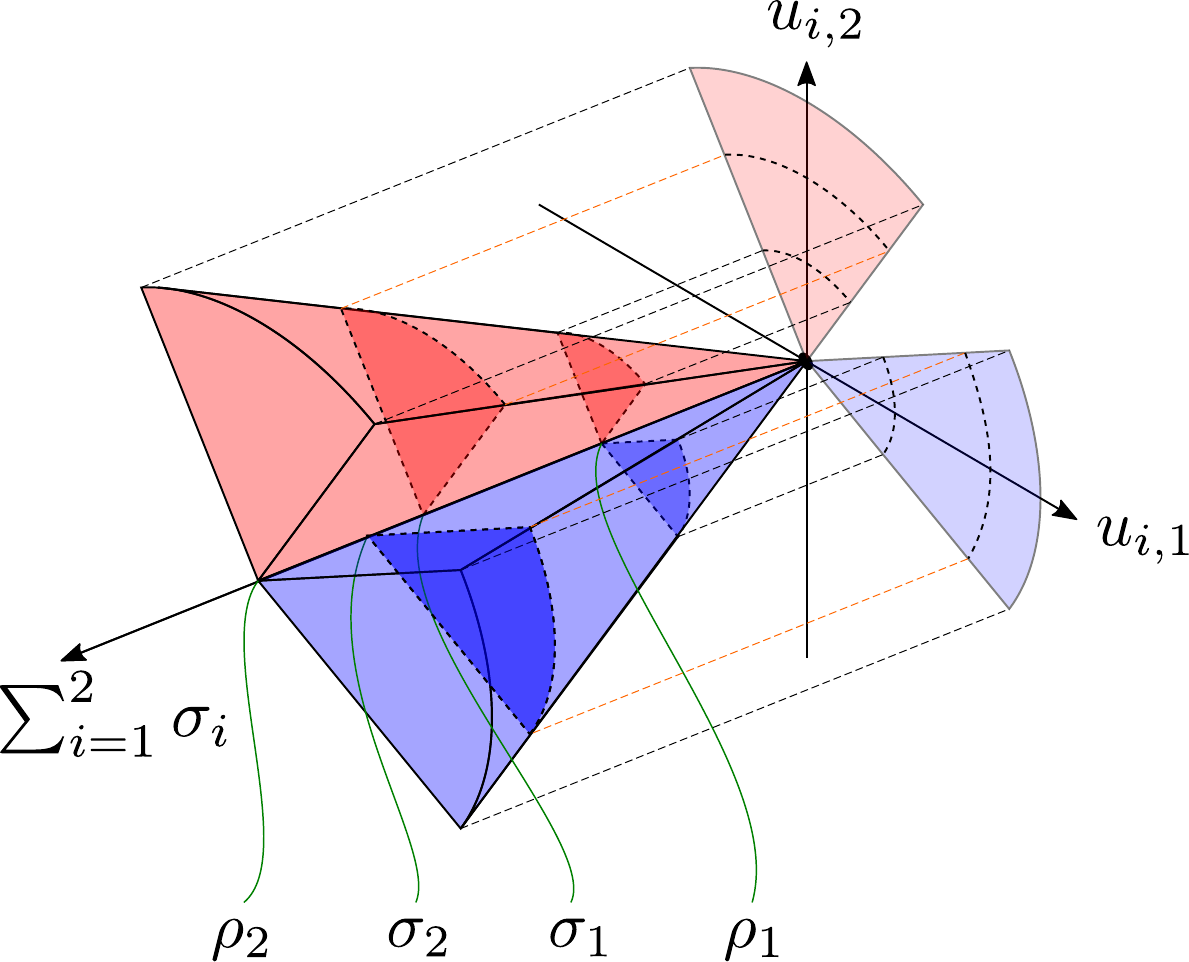}
    \caption{The non-convexity is removed by relaxing \eqref{eq:ocp_c} to
      \eqref{eq:rcp_c}-\eqref{eq:rcp_d}.}
    \label{fig:input_set_relaxation_2}
  \end{subfigure}%
  \hfill%
  \begin{subfigure}[b]{0.32\textwidth}
    \centering
    \includegraphics[width=0.9\columnwidth]{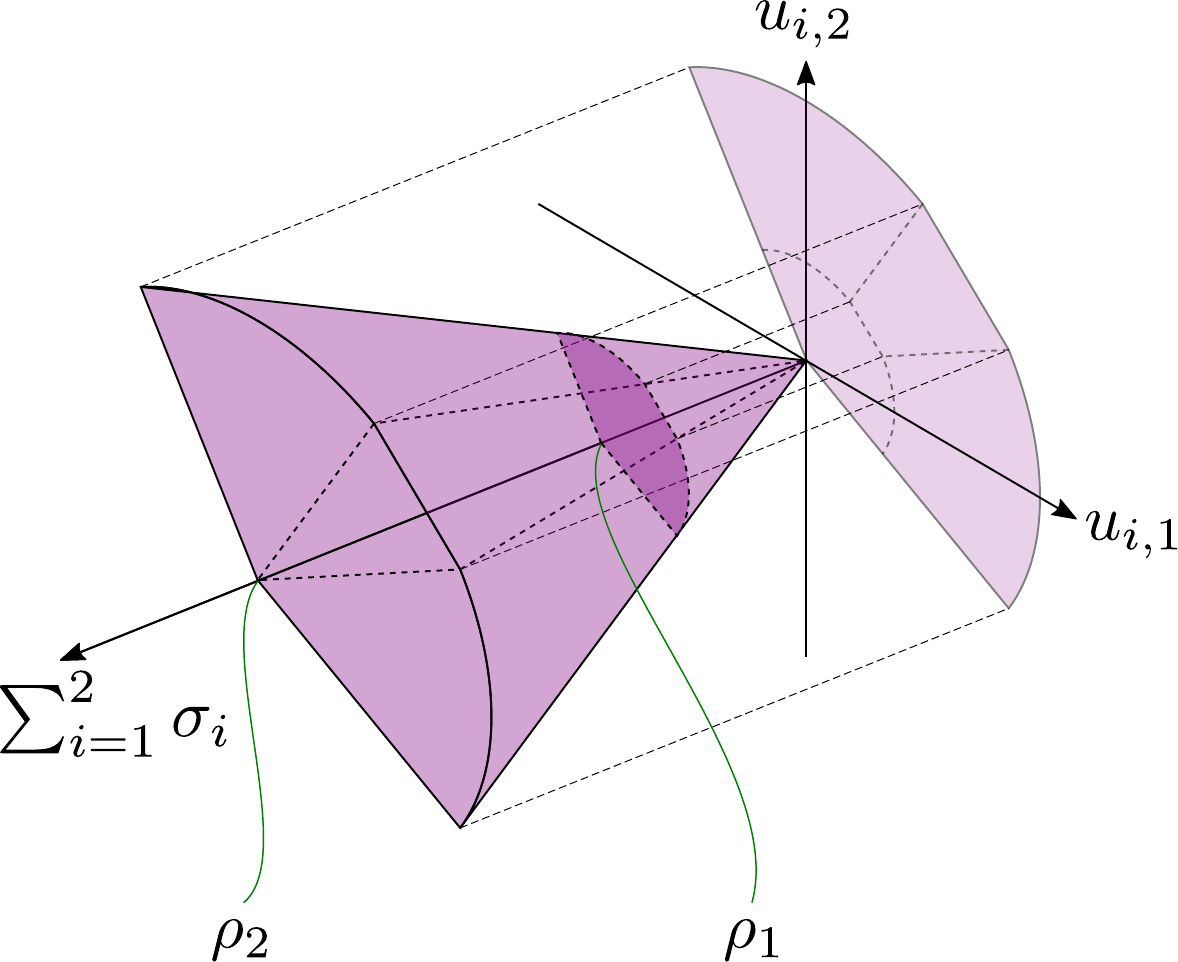}
    \caption{The mutual exclusivity is removed by relaxing \eqref{eq:ocp_d} to
      \eqref{eq:rcp_e}.}
    \label{fig:input_set_relaxation_3}
  \end{subfigure}
  \caption{Problem~\ref{problem:rcp} convexifies the input set of
    Problem~\ref{problem:ocp}, here shown for $M=2$, $K=1$ and $m=2$. The
    relaxation consists of three steps:
    \protect\subref{fig:input_set_relaxation_1})
    \eqref{eq:ocp_c}-\eqref{eq:ocp_f} originally constrain each input to
    non-convex disjoint sets;
    \protect\subref{fig:input_set_relaxation_2}) by relaxing \eqref{eq:ocp_c} to
    \eqref{eq:rcp_c}-\eqref{eq:rcp_d}, individual input sets are convexified to
    3D slices; \protect\subref{fig:input_set_relaxation_3}) by relaxing
    \eqref{eq:ocp_d} to \eqref{eq:rcp_e}, a convex hull is obtained.}
  \label{fig:input_set_relaxation}
\end{figure*}

This section presents Theorem~\ref{theorem:lcvx}, which is the main result of
the paper and states that Problem~\ref{problem:ocp} can be solved via convex
optimization under certain conditions.

Problem~\ref{problem:ocp} is non-convex due to the input norm lower-bound in
\eqref{eq:ocp_c} and is mixed-integer due to
\eqref{eq:ocp_d}. Figure~\ref{fig:input_set_relaxation_1} illustrates the
non-convex mixed-integer nature of the problem due to the input sets being
non-convex and disjoint. Consider the following convex relaxation of
Problem~\ref{problem:ocp}:
\begin{poptimization}{left}{rcp}{}{u_i,\gamma_i,\sigma_i,t_f}{%
    m(t_f,x(t_f))
  }{$\mathcal R$}
  \dot x(t) = Ax(t)+B\inlinesum_{i=1}^M{u_i(t)}+w,~x(0) = x_0, \#
  \gamma_i(t)\rho_1 \le \sigma_i(t)\le \gamma_i(t)\rho_2\quad i=1,\dots,M, \#
  \|u_i(t)\|_2\le\sigma_i(t)\quad i=1,\dots,M, \#
  0\le \gamma_i(t)\le 1\quad i=1,\dots,M, \#
  \inlinesum_{i=1}^M\gamma_i(t) \le K, \#
  C_iu_i(t)\le 0\quad i=1,\dots,M, \#
  b(t_f,x(t_f)) = 0.
\end{poptimization}

\begin{figure}
  \centering
  \includegraphics[width=1\columnwidth]{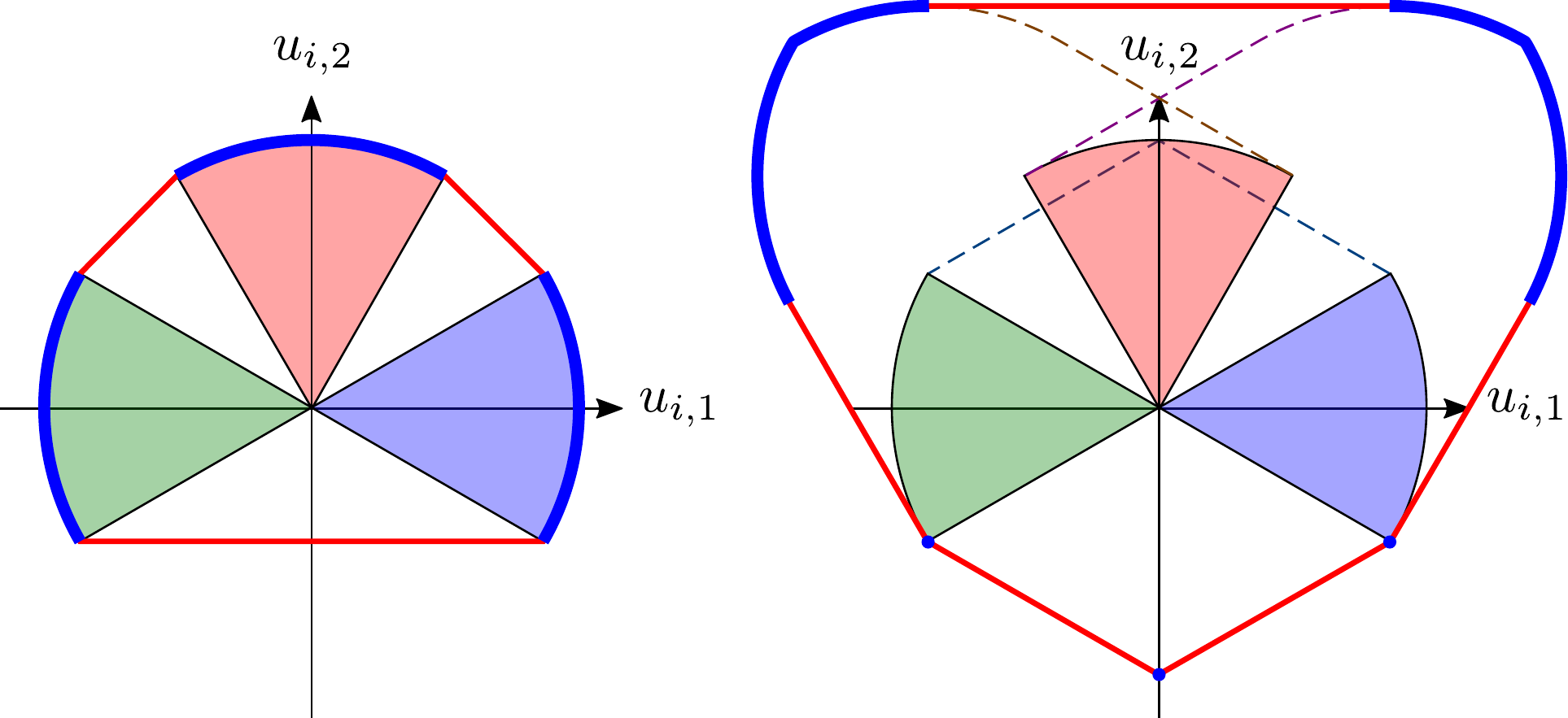}
  \caption{The relaxed input set of Problem~\ref{problem:rcp} is the convex hull
    of the Minkowski sums for every combination of $K$ or fewer of the
    individual input sets (each one relaxed via
    \eqref{eq:rcp_c}-\eqref{eq:rcp_e}). For $M=3$, $K=2$ and $m=2$, (left) shows
    the case of one- and (right) shows the case of two-input set combinations,
    with the relaxed set shown in bold red and the boundaries of the constituent
    Minkowski sums shown as dashed lines. The optimal solution takes values from
    the extreme points, shown as blue segments. The origin is also an extreme
    point, corresponding to a combination of zero input sets (not shown).}
  \label{fig:minkowski_sum}
\end{figure}

By replacing \eqref{eq:ocp_c}-\eqref{eq:ocp_d} with
\eqref{eq:rcp_c}-\eqref{eq:rcp_e}, the input set of Problem~\ref{problem:rcp}
becomes the convex hull of the Minkowski sums of every combination of at most
$K$ of the relaxed individual input sets of
Problem~\ref{problem:ocp}. Figure~\ref{fig:input_set_relaxation} illustrates the
convex relaxation for $M=2$ and $K=1$, where it is possible to also visualize
the $\sigma_i$ slack variables to elucidate the convex
lifting. Figure~\ref{fig:minkowski_sum} shows the case of $M=3$ and $K=2$, where
the input set is projected onto the $u_i$ space. It will be shown in
Section~\ref{sec:ct_lcvx_proof} that the optimal solution is extremal, hence it
will take values among the extreme points of the input set of
Problem~\ref{problem:rcp} with at most $K$ inputs active.

Consider the following conditions, which are sufficient to eliminate degenerate
optimal solutions of Problem~\ref{problem:rcp} that may be infeasible for
Problem~\ref{problem:ocp}. To state the conditions, define an \textit{adjoint
  system} whose output $y(t)\in\reals^m$ is called the \textit{primer vector}:
\begin{eqnarray}
  \label{eq:adjoint_system}
  \dot\lambda(t)=-A\T\lambda(t),\quad y(t)=B\T\lambda(t).
\end{eqnarray}

Furthermore, it will be seen in the proof of Lemma~\ref{lemma:lcvx} that we are
interested in ``how much'' $y(t)$ projects onto the $i$-th input pointing
set. This is given by the following input \textit{gain} measure:
\begin{eqnarray}
  \label{eq:input_gain}
  \Gamma_i(t) \definedas \Proj{\set{U_i}}{y(t)}.
\end{eqnarray}

\begin{condition}{}
  \label{condition:observability}
  The adjoint system \eqref{eq:adjoint_system} is observable.
\end{condition}

\begin{condition}{}
  \label{condition:normality}
  The adjoint system \eqref{eq:adjoint_system} and pointing cone geometry
  \eqref{eq:ocp_f} satisfy either:
  \begin{enumerate}
  \item[(a)] $\Gamma_i(t)\ne 0~\alev{[0,t_f]}$ $\forall i$
    s.t. $y(t)\notin\interior{\normalcone{\set U_i}{0}}$;
  \item[(b)] on any interval where $\Gamma_i(t)=0$, $\Gamma_j(t)>0$ for at least
    $K$ other inputs.
  \end{enumerate}
\end{condition}

\begin{condition}{}
  \label{condition:ambiguity}
  The adjoint system \eqref{eq:adjoint_system} and pointing cone geometry
  \eqref{eq:ocp_f} satisfy either:
  \begin{enumerate}
  \item[(a)] $\Gamma_i(t)\ne \Gamma_j(t)~\alev{[0,t_f]}$ $\forall i$
    s.t. $y(t)\notin\interior{\normalcone{\set U_i}{0}}$;
  \item[(b)] on any interval where $\Gamma_i(t)=\Gamma_j(t)$, there exist $K$
    inputs with $\Gamma_k(t)>\Gamma_i(t)$ or $M-K$ inputs where
    $\Gamma_k(t)<\Gamma_i(t)$.
  \end{enumerate}
\end{condition}

\begin{condition}{}
  \label{condition:transversality}
  The following intersection holds:
  \begin{eqnarray}
    \label{eq:transversality_intersection}
    \range
    \begin{bmatrix}
      \grad_x b[t_f]\T \\
      \grad_t b[t_f]\T
    \end{bmatrix}
    \intersection
    \conehull
    \begin{bmatrix}
      \grad_x m[t_f]\T \\
      \grad_t m[t_f]
    \end{bmatrix} = \{0\}.
  \end{eqnarray}
\end{condition}

We now state the main result of this paper, which claims that
Problem~\ref{problem:rcp} solves Problem~\ref{problem:ocp} under the above
conditions. The theorem is proved in Section~\ref{sec:ct_lcvx_proof}.

\begin{theorem}{}
  \label{theorem:lcvx}
  The solution of Problem~\ref{problem:rcp} is globally optimal $\alev{[0,t_f]}$
  for Problem~\ref{problem:ocp} if
  Conditions~\ref{condition:observability}-\ref{condition:transversality} hold.
\end{theorem}

\subsection{Discussion of
  Conditions~\ref{condition:observability}-\ref{condition:transversality} and
  Special Cases}
\label{subsec:condition_special_cases}

\begin{figure}
  \centering
  \begin{subfigure}[b]{\columnwidth}
    \centering
    \includegraphics[width=0.7\textwidth]{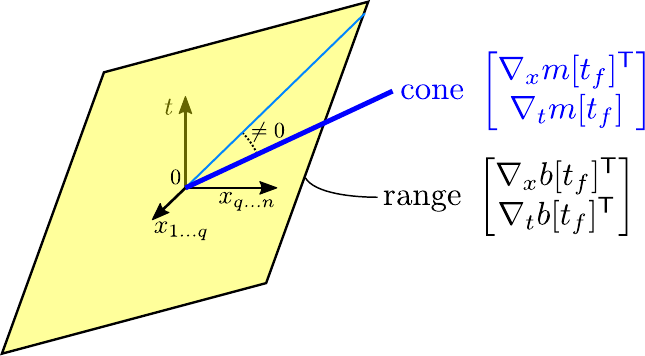}
    \caption{General Condition~\ref{condition:transversality} requirement.}
    \label{fig:condition_transversality}
  \end{subfigure}%

  \begin{subfigure}[b]{0.48\columnwidth}
    \centering
    \includegraphics[width=0.8\columnwidth]{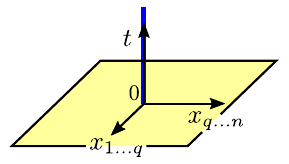}
    \caption{Minimum-time.}
    \label{fig:condition_transversality_case_1}
  \end{subfigure}%
  \hfill%
  \begin{subfigure}[b]{0.48\columnwidth}
    \centering
    \includegraphics[width=0.8\columnwidth]{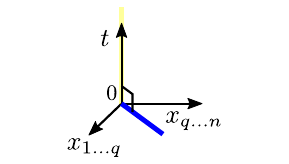}
    \caption{Minimum-error.}
    \label{fig:condition_transversality_case_2}
  \end{subfigure}
  \caption{Condition~\ref{condition:transversality} for the general and two
    special cases.}
\end{figure}

This section describes situations when
Conditions~\ref{condition:observability}-\ref{condition:transversality} are easy
to verify. First, Condition~\ref{condition:observability} is readily verified by
checking if the pair $\{-A\T,B\T\}$ in \eqref{eq:adjoint_system} is observable
\cite{Antsaklis2007}.

Condition~\ref{condition:transversality} originates from the maximum principle
transversality requirement
\eqref{eq:ocp_nonsmooth_transversality}. Figure~\ref{fig:condition_transversality}
illustrates the general requirement. As shown in
Figure~\ref{fig:condition_transversality_case_1}, the condition is satisfied for
a minimum-time problem with a fixed or (partially) free terminal state. As shown
in Figure~\ref{fig:condition_transversality_case_2}, the condition is also
satisfied for a fixed- or free-time problem with a (partially) penalized
terminal state. Note that Assumption~\ref{ass:full_range_cone} must hold.

\begin{figure*}
  \centering
  \begin{subfigure}[b]{1\columnwidth}
    \centering
    \includegraphics[width=0.48\textwidth]{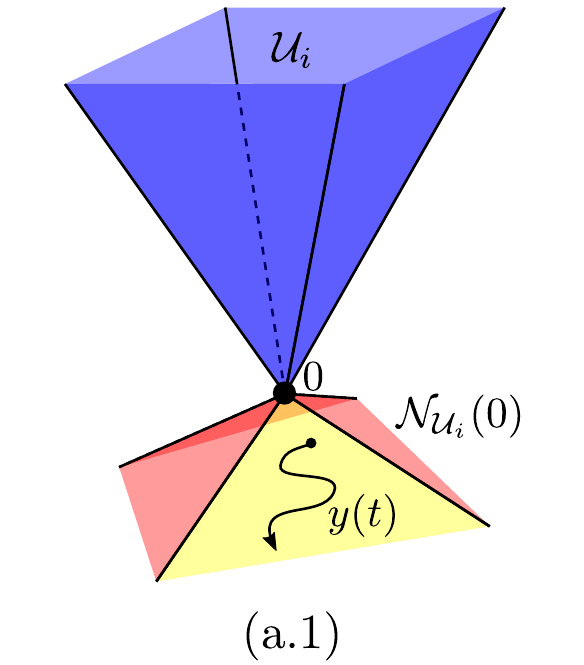}%
    \hfill%
    \includegraphics[width=0.48\textwidth]{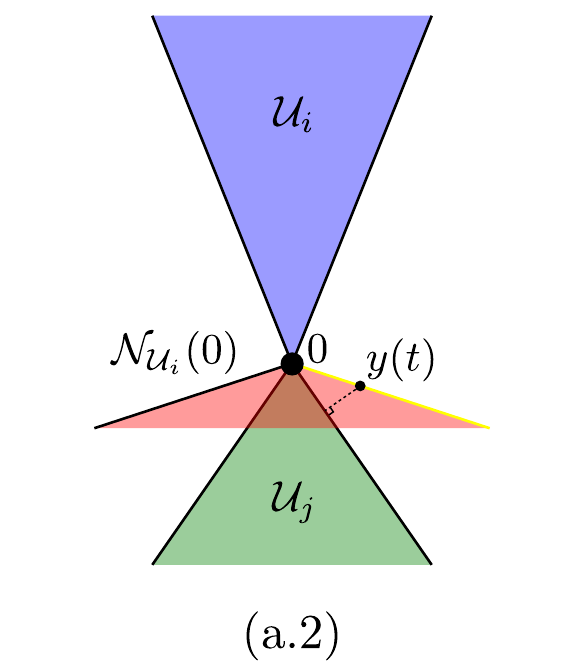}
    \caption{When case (a) of Condition~\ref{condition:normality} fails, $y(t)$
      can evolve on the normal cone boundary (a.1). Case (b) then holds if
      $y(t)\in\boundary{\normalcone{\mathcal U_i}{0}}$ projects positively onto
      at least $K$ other input pointing sets. For $K=1$, (a.2) illustrates a
      case where $y(t)$ projects positively onto $\mathcal U_j$.}
    \label{fig:condition_facet_degeneracy}
  \end{subfigure}%
  \hfill%
  \begin{subfigure}[b]{1\columnwidth}
    \centering
    \includegraphics[width=0.48\textwidth]{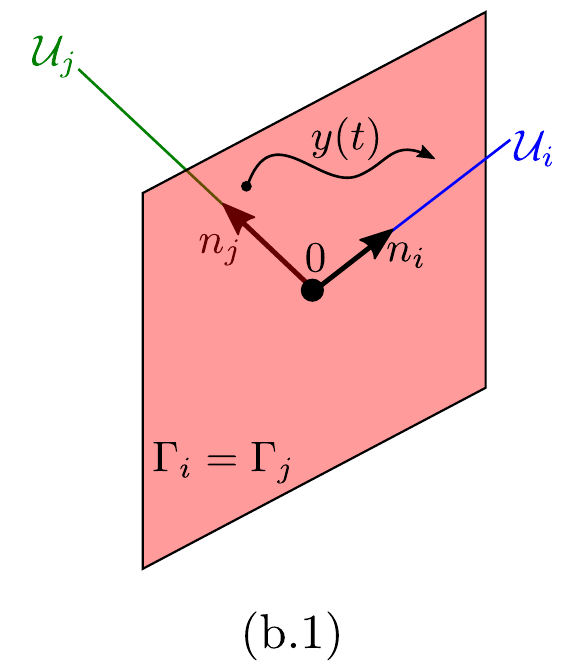}%
    \hfill%
    \includegraphics[width=0.48\textwidth]{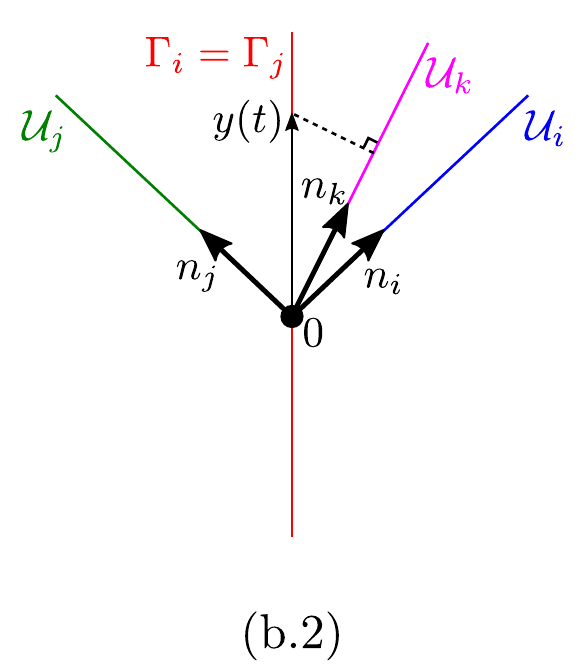}
    \caption{When case (a) of Condition~\ref{condition:ambiguity} fails, $y(t)$
      can evolve on the manifold $\Gamma_i(t)=\Gamma_j(t)$ as shown in (b.1) for
      ray cones. Case (b) then holds if, for example, $y(t)$ projects more
      positively onto $K$ other cones. For $K=1$, (b.2) shows a situation where
      $y(t)$ projects more positively onto $\mathcal U_k$.}
    \label{fig:condition_projection_ambiguity}
  \end{subfigure}
  \caption{Visualization of Conditions~\ref{condition:normality} and
    \ref{condition:ambiguity} when their case (a) fails. In such circumstances,
    the conditions can nevertheless hold given the right input pointing set
    geometry.}
  \label{fig:condition_illustrations}
\end{figure*}

Conditions~\ref{condition:normality} and \ref{condition:ambiguity} are
illustrated in Figure~\ref{fig:condition_illustrations}. Both conditions can be
checked via matrix algebra in the special case of \textit{ray} cones
$\set U_i=\conehull n_i$ for some \textit{direction} $n_i\in\reals^m$. In this
case, $\boundary{\normalcone{\set U_i}{0}}=\{u\in\reals^m:n_i\T u=0\}$ which is
a hyperplane. Consider the adjoint system \eqref{eq:adjoint_system} with the
``projected'' primer vector:
\begin{eqnarray}
  \label{eq:projected_primer_vector}
  n_i\T y(t)=(Bn_i)\T\lambda(t).
\end{eqnarray}

If the pair $\{-A\T,(Bn_i)\T\}$ is observable then
Condition~\ref{condition:normality} case (a) holds. If not, let $\set V_i$ be
the unobservable subspace and consider the special case
$\range B\T(-A\T)^k\set V_i\subseteq\range z_i$ $\forall k=0,\dots,n-1$ for some
$z_i\in\real^m$. If $\Proj{\set U_k}{z_i}>0$ and $\Proj{\set U_k}{-z_i}>0$ for
$K$ or more input pointing sets, then Condition~\ref{condition:normality} case
(b) holds. Similarly, for Condition~\ref{condition:ambiguity} we consider the
projected primer vector:
\begin{eqnarray}
  \label{eq:projected2_primer_vector}
  (n_i-n_j)\T y(t)=\big(B(n_i-n_j)\big)\T\lambda(t).
\end{eqnarray}

If the pair $\{-A\T,\big(B(n_i-n_j)\big)\T\}$ is observable then
Condition~\ref{condition:ambiguity} case (a) holds. If not, let $z_i\in\reals^m$
be defined as before with $\set V_i$ the unobservable subspace for this new
system. If $\Proj{\set U_k}{z_i}>\Proj{\set U_i}{z_i}$ for $K$ other inputs or
$\Proj{\set U_k}{z_i}<\Proj{\set U_i}{z_i}$ for $M-K$ other inputs, then
Condition~\ref{condition:ambiguity} case (b) holds.

\section{Nonsmooth Maximum Principle}
\label{sec:pmp}

This section states a nonsmooth version of the maximum principle that we shall
use for proving Theorem~\ref{theorem:lcvx}. Consider the following general
optimal control problem:
\begin{poptimization}{left}{ocp_general}{}{u,t_f}{%
    m(t_f,x(t_f))}{$\mathcal G$}
  \dot x(t) = f(t,x(t),u(t)),\quad x(\tzero)=x_0, \#
  g(t,u(t))\le 0, \#
  b(t_f,x(t_f)) = 0.
\end{poptimization}
where the state trajectory $x(\cdot)$ is absolutely continuous and the control
trajectory $u(\cdot)$ is measurable. The dynamics
$f:\reals\times\reals^n\times\reals^m\to\reals^n$ are convex and continuously
differentiable. The terminal cost $m:\reals\times\reals^n\to\reals$, the input
constraint $g:\reals\times\reals^m\to\reals^{n_g}$, and the terminal constraint
$b:\reals\times\reals^n\to\reals^{n_b}$ are convex. Let us denote the terminal
manifold as
$\set T\definedas\{x\in\reals^n:\textnormal{\eqref{eq:ocp_general_d} holds}\}$. The
Hamiltonian function is defined as:
\begin{eqnarray}
  \label{eq:ocp_nonsmooth_hamiltonian}
  H(t,x(t),u(t),\lambda(t))\definedas \lambda(t)\T f[t],
\end{eqnarray}
where $\lambda(\cdot)$ is the \textit{adjoint variable} trajectory. We now state
the nonsmooth maximum principle, due to \cite[Theorem~8.7.1]{Vinter2000} (see
also \cite{Clarke2010,Hartl1995}), which specifies the necessary conditions of
optimality for Problem~\ref{problem:ocp_general}.

\begin{theorem}{}[Maximum Principle]
  \label{theorem:mp}
  Let $x(\cdot)$ and $u(\cdot)$ be optimal on the interval $[0,t_f]$. There
  exist an \textit{abnormal multiplier} $\alpha\le 0$ and an absolutely
  continuous $\lambda(\cdot)$ such that the following conditions are satisfied:
  \begin{enumerate}
  \item Non-triviality:
    \begin{eqnarray}
      \label{eq:ocp_nonsmooth_non_triviality}
      (\alpha,\lambda(t))\ne 0~\forall t\in [\tzero,t_f];
    \end{eqnarray}
  \item Pointwise maximum:
    \begin{eqnarray}
      \label{eq:ocp_nonsmooth_pointwise_maximum_argmax}
      u(t) = \argmax_{v\in\textnormal{\eqref{eq:ocp_general_c}}}H(t,x(t),v,\lambda(t))~\alev{[0,t_f]};
    \end{eqnarray}
  \item The differential equations:
    \begin{subequations}
      \begin{flalign}
        \label{eq:ocp_nonsmooth_dynamics_state}
        &\dot x(t)
        = \grad_{\lambda} H[t]\T~\alev{[\tzero,t_f]}, &&\\
        \label{eq:ocp_nonsmooth_dynamics_adjoint}
        &\dot \lambda(t)
        = -\grad_x H[t]\T~\alev{[\tzero,t_f]}, &&\\
        \label{eq:ocp_nonsmooth_dynamics_hamiltonian}
        &\dot H[t]
        = \grad_t H[t]~\alev{[\tzero,t_f]};
      \end{flalign}
    \end{subequations}
  \item Transversality:
    \begin{subequations}
      \label{eq:ocp_nonsmooth_transversality}
      \begin{flalign}
        \label{eq:ocp_nonsmooth_transversality_adjoint}
        &\lambda(t_f) \in \alpha\subdiff[x] m[t_f]\T+\normalcone{\set T}{x(t_f)}, &&\\
        \label{eq:ocp_nonsmooth_transversality_hamiltonian}
        &0 \in H[t_f]+\alpha\subdiff[t] m[t_f]+\normalcone{\set T}{t_f}.
      \end{flalign}
    \end{subequations}
  \end{enumerate}
\end{theorem}

\section{Lossless Convexification Proof}
\label{sec:ct_lcvx_proof}

This section proves Theorem~\ref{theorem:lcvx} in two steps. First, the output
of Problem~\ref{problem:rcp} is shown to be feasible for
Problem~\ref{problem:ocp} via a maximum principle argument. Second, this
solution is shown to also be globally optimal via an equal cost function
argument.

\begin{lemma}
  \label{lemma:lcvx}
  The solution of Problem~\ref{problem:rcp} is feasible $\alev{[0,t_f]}$ for
  Problem~\ref{problem:ocp} if Conditions
  \ref{condition:observability}-\ref{condition:transversality} hold.
\end{lemma}

\begin{pf*}{Proof.}
  The proof uses the maximum principle from Theorem~\ref{theorem:mp}. For
  Problem~\ref{problem:rcp}, the adjoint and Hamiltonian dynamics follow from
  \eqref{eq:ocp_nonsmooth_dynamics_adjoint} and
  \eqref{eq:ocp_nonsmooth_dynamics_hamiltonian}:
  \begin{subequations}
    \label{eq:proof_dynamics}
    \begin{flalign}
      &\dot\lambda(t) = -A\T\lambda(t)
      ~\alev{[0,t_f]},\hspace{-5mm} \label{eq:lambda_dynamics} &&\\
      &\dot H[t] = 0~\alev{[0,t_f]}, \label{eq:hamiltonian_dynamics}
    \end{flalign}
  \end{subequations}

  Using the subdifferential basic chain rule
  \cite[Theorem~10.6]{Rockafellar1998}, the transversality condition
  \eqref{eq:ocp_nonsmooth_transversality} yields:
  \begin{subequations}
    \label{eq:proof_transversality}
    \begin{flalign}
      \label{eq:proof_transversality_lambda}
      &\lambda(t_f) = \grad_x m[t_f]\T\alpha+\grad_xb[t_f]\T\beta, &&\\
      \label{eq:proof_transversality_H}
      &H[t_f] = -\grad_t m[t_f]\alpha-\grad_t b[t_f]\T\beta,
    \end{flalign}
  \end{subequations}
  for some $\beta\in\reals^{n_b}$. Due to \eqref{eq:hamiltonian_dynamics},
  \eqref{eq:proof_transversality_H} and absolute continuity, we have
  \cite[Theorem~9]{Varberg1965}:
  \begin{eqnarray}
    \label{eq:hamiltonian_constant}
    H[t]=-\grad_t m[t_f]\alpha-\grad_t b[t_f]\T\beta,~\forall t\in[0,t_f].
  \end{eqnarray}

  We claim that the primer vector $y(t)\ne 0~\alev{[0,t_f]}$. Since $y(t)$ is
  the output of \eqref{eq:adjoint_system}, it is an analytic function and
  $y(t)=0$ either $\forall t\in [0,t_f]$ or at a countable number of instances
  \cite{Acikmese2011,Carson2011}. By contradiction, suppose that $y(t)=0$
  $\forall t\in [0,t_f]$. Since Condition~\ref{condition:observability} holds,
  $\lambda(0)=0$. Since \eqref{eq:lambda_dynamics} is homogeneous,
  $\lambda(t)=0$ $\forall t\in[0,t_f]$. The transversality condition
  \eqref{eq:proof_transversality} hence simplifies to:
  \begin{eqnarray}
    \label{eq:transversality_y0}
    \begin{bmatrix}
      \grad_x m[t_f]\T \\
      \grad_t m[t_f]
    \end{bmatrix}(-\alpha) =
    \begin{bmatrix}
      \grad_x b[t_f]\T \\
      \grad_t b[t_f]\T
    \end{bmatrix}\beta,
  \end{eqnarray}
  and since Condition~\ref{condition:transversality} holds, $\alpha=0$. Hence
  $(\alpha,\lambda(t))=0$ $\forall t\in [0,t_f]$, which violates
  non-triviality \eqref{eq:ocp_nonsmooth_non_triviality}. Therefore it must be
  that $y(t)\ne 0~\alev{[0,t_f]}$. Because the necessary conditions are
  scale-invariant, we can set $\alpha=-1$ without loss of generality. The
  pointwise maximum condition
  \eqref{eq:ocp_nonsmooth_pointwise_maximum_argmax} implies that the following
  must hold $\alev{[0,t_f]}$: \def\opticmd{\argmax}
  \begin{optimization}{left}{argmax_1}{}{u_i,\gamma_i,\sigma_i}{%
      \begin{array}{l}
        \inlinesum_{i=1}^My(t)\T u_i(t)
      \end{array}}{}
    \textnormal{constraints \eqref{eq:rcp_c}-\eqref{eq:rcp_g} hold.}
  \end{optimization}

  We shall now analyze the optimality conditions of \eqref{eq:argmax_1}. For
  concise notation, the time argument $t$ shall be omitted. Expressing
  \eqref{eq:argmax_1} as a minimization and treating constraints
  \eqref{eq:rcp_e} and \eqref{eq:rcp_f} implicitly, we can write the
  Lagrangian of \eqref{eq:argmax_1} \cite{Boyd2004}:
  \begin{flalign}
    \label{eq:lagrangian}
    &\mathcal L(u_i,\gamma_i,\sigma_i,\lambda_{1\dots 4}^i) = \inlinesum_{i=1}^M
    -y\T u_i+\lambda_1^i(\norm{u_i}{2}-\hspace{-5mm} &&\\
    &\hspace{10mm}\sigma_i)+\lambda_2^i(\gamma_i\rho_1-\sigma_i)+
    \lambda_3^i(\sigma_i-\gamma_i\rho_2)+{\lambda_4^i}\T C_iu_i, \nonumber
  \end{flalign}
  where $\lambda_j^i\ge 0$ are Lagrange multipliers satisfying the following
  complementarity conditions:
  \begin{subequations}
    \label{eq:complementarity_conditions}
    \begin{flalign}
      &\lambda_1^i(\norm{u_i}{2}-\sigma_i) = 0, \label{eq:kkt_1} &&\\
      &\lambda_2^i(\gamma_i\rho_1-\sigma_i) = 0, \label{eq:kkt_2} &&\\
      &\lambda_3^i(\sigma_i-\gamma_i\rho_2) = 0, \label{eq:kkt_3} &&\\
      &\lambda_4^i\circ C_iu_i = 0. \label{eq:kkt_4}
    \end{flalign}
  \end{subequations}

  Next, the Lagrange dual function is given by:
  \begin{flalign}
    &g(\lambda_{1\dots 4}^i) = \inf_{u_i,\gamma_i,\sigma_i}
    \mathcal L(u_i,\gamma_i,\sigma_i,\lambda_{1\dots 4}^i) \nonumber &&\\
    &\phantom{g(\lambda_{1\dots 4}^i) } =\inlinesum_{i=1}^M\inf_{\sigma_i}\left[
      (\lambda_3^i-\lambda_2^i-\lambda_1^i)\sigma_i\right]- \nonumber &&\\
    &\phantom{g(\lambda_{1\dots 4}^i) =}\inlinesum_{i=1}^M\sup_{u_i}\left[(y-C_i\T\lambda_4^i)\T u_i-
      \lambda_1^i\norm{u_i}{2}\right]+ \nonumber &&\\
    &\phantom{g(\lambda_{1\dots 4}^i) =}\textstyle\inf_{\textnormal{\eqref{eq:rcp_e},\eqref{eq:rcp_f}}}\inlinesum_{i=1}^M
    (\lambda_2^i\rho_1-\lambda_3^i\rho_2)\gamma_i.
    \label{eq:dual_function}
  \end{flalign}
  

  The dual function bounds the primal optimal cost from above. A non-trivial
  upper-bound requires:
  \begin{subequations}
    \begin{flalign}
      &\norm{y-C_i\T\lambda_4^i}{2}\le\lambda_1^i, \label{eq:dual_function_result_1} &&\\
      &\lambda_3^i-\lambda_2^i-\lambda_1^i = 0, \label{eq:dual_function_result_2}
    \end{flalign}
  \end{subequations}
  where the first inequality is akin to the $\norm{\cdot}{2}$ conjugate
  function \cite[Example~3.26]{Boyd2004}. However, note that if
  \eqref{eq:dual_function_result_1} is strict then $\norm{u_i}{2}=0$ is
  optimal, which is trivially feasible for
  Problem~\ref{problem:ocp}. Substituting \eqref{eq:dual_function_result_2}
  into \eqref{eq:dual_function_result_1} gives the following condition for
  non-trivial solutions:
  \begin{eqnarray}
    \label{eq:prelim_characteristic_equation_1}
    \norm{y-C_i\T\lambda_4^i}{2} = \lambda_3^i-\lambda_2^i.
  \end{eqnarray}

  Next, note that a non-trivial solution implies $\gamma_i>0$. Due to
  Assumption~\ref{ass:thrust_bounds}, \eqref{eq:kkt_2} and \eqref{eq:kkt_3}, a
  non-trivial solution cannot have $\lambda_2^i>0$ and $\lambda_3^i>0$
  simultaneously. Furthemore, \eqref{eq:dual_function} reveals that
  $\gamma_i>0$ is not sub-optimal if and only if
  $\lambda_2^i\rho_1-\lambda_3^i\rho_2\le 0$. Hence $\lambda_2^i=0$ and
  $\lambda_3^i\ge 0$ are necessary for optimality. As a result
  \eqref{eq:prelim_characteristic_equation_1} simplifies to:
  \begin{eqnarray}
    \label{eq:prelim_characteristic_equation_2}
    \norm{y-C_i\T\lambda_4^i}{2} = \lambda_3^i.
  \end{eqnarray}

  Next, note that at optimality the left-hand side of
  \eqref{eq:prelim_characteristic_equation_2} equals the Eucledian projection
  onto $\set{U_i}$,
  i.e. $\norm{y-C_i\T\lambda_4^i}{2}=\Proj{\set{U_i}}{y}$. This can be shown
  by contradiction using Assumption~\ref{ass:full_range_cone},
  \eqref{eq:kkt_4} and that it is optimal to choose
  $u_i=\norm{u_i}{2}(y-C_i\T\lambda_4^i)/\norm{y-C_i\T\lambda_4^i}{2}$ in
  \eqref{eq:dual_function}. Note that the degenerate case of $u_i\ne 0$ and
  $\norm{y-C_i\T\lambda_4^i}{2}=0$ is eliminated by
  Condition~\ref{condition:normality}, as discussed below. Thus
  \eqref{eq:prelim_characteristic_equation_2} simplifies to the following
  relationship, which we call the \textit{characteristic equation} of
  non-trivial solutions to \eqref{eq:argmax_1}:
  \begin{eqnarray}
    \label{eq:characteristic_equation}
    \Proj{\set{U_i}}{y} = \lambda_3^i.
  \end{eqnarray}

  Note that when $\lambda_3^i>0$ then $\norm{u_i}{2}=\sigma_i=\gamma_i\rho_2$
  due to \eqref{eq:kkt_3} and \eqref{eq:dual_function_result_2}. Substituting
  \eqref{eq:characteristic_equation} into \eqref{eq:dual_function} yields:
  \begin{eqnarray}
    \label{eq:dual_function_simplified}
    g(\lambda_{1\dots 4}^i) =
    -\rho_2\textstyle\sup_{\textnormal{\eqref{eq:rcp_e},\eqref{eq:rcp_f}}}\inlinesum_{i=1}^{K'}
    \Proj{\set{U_i}}{y}\gamma_i,
  \end{eqnarray}
  where we assume that the characteristic equation
  \eqref{eq:characteristic_equation} does not hold for $i=K'+1,\dots,M$ such
  that $\gamma_{i>K'}=0$. To facilitate discussion, define the $i$-th input
  \textit{gain} as in \eqref{eq:input_gain}. Note that $\Gamma_i\ge 0$ due to
  \eqref{eq:characteristic_equation}. Thus \eqref{eq:dual_function_simplified}
  becomes:
  \begin{eqnarray}
    \label{eq:dual_function_gain}
    g(\lambda_{1\dots 4}^i) =
    -\rho_2\textstyle\sup_{\textnormal{\eqref{eq:rcp_e},\eqref{eq:rcp_f}}}
    \inlinesum_{i=1}^{K'}\Gamma_i\gamma_i.
  \end{eqnarray}

  Without loss of generality, assume a descending ordering
  $\Gamma_i\ge\Gamma_j$ for $i>j$. Let $K''\definedas\min\{K,K'\}$. By
  inspection of \eqref{eq:dual_function_gain}, the condition:
  \begin{eqnarray}
    \label{eq:feasibility_condition}
    \Gamma_{K''}>0~\land~\Gamma_{K''}>\Gamma_{K''+1},
  \end{eqnarray}
  is sufficient to ensure that it is optimal to set
  \begin{eqnarray}
    \label{eq:gamma_optimality_structure}
    \gamma_i=
    \begin{cases}
      1 & \text{for }i\le K'', \\
      0 & \text{otherwise.}
    \end{cases}
  \end{eqnarray}

  The lemma is proved if \eqref{eq:feasibility_condition} holds
  $\alev{[0,t_f]}$. However, this holds by assumption due to
  Conditions~\ref{condition:normality} and \ref{condition:ambiguity}. In
  particular, Condition~\ref{condition:normality} case (a) assures that
  $\Gamma_{K''}>0~\alev{[0,t_f]}$. If on some interval $\Gamma_k=0$,
  Condition~\ref{condition:normality} case (b) assures that $k>K''$. Next, if
  $K''<K$ then due to $\Gamma_{K''}>0$ and the definition of $K'$, it must be
  that $\Gamma_{K''+1}=0\Rightarrow \Gamma_{K''}>\Gamma_{K''+1}$. On the other
  hand, if $K''=K$ then Condition~\ref{condition:ambiguity} case (a) assures
  that $\Gamma_{K}>\Gamma_{K+1}~\alev{[0,t_f]}$. If on some interval
  $\Gamma_k=\Gamma_{k+1}$, Condition~\ref{condition:ambiguity} case (b) assures
  that $k\ne K$.

  


  Thus, \eqref{eq:feasibility_condition} holds $\alev{[0,t_f]}$ and the lemma
  is proved. From \eqref{eq:gamma_optimality_structure}, the structure of the
  optimal solution is \textbf{bang-bang with at most $K$ inputs active}
  $\alev{[0,t_f]}$. \qed
\end{pf*}

Lemma~\ref{lemma:lcvx} guarantees that solving Problem~\ref{problem:rcp} yields
a feasible solution of Problem~\ref{problem:ocp}. We will now show that this
solution is globally optimal, thus proving Theorem~\ref{theorem:lcvx}.

\begin{pf*}{Proof of Theorem~\ref{theorem:lcvx}.}
  The solution of Problem~\ref{problem:rcp} is feasible $\alev{[0,t_f]}$ for
  Problem~\ref{problem:ocp} due to Lemma~\ref{lemma:lcvx}. Furthermore, the cost
  functions of Problems~\ref{problem:ocp} and \ref{problem:rcp} are the
  same. The optimal costs thus satisfy
  $\optimal{J_{\mathcal O}}\le\optimal{J_{\mathcal R}}$. However, any solution
  of Problem~\ref{problem:ocp} is feasible for Problem~\ref{problem:rcp} by
  setting $\sigma_i(t)=\norm{u_i(t)}{2}$, thus
  $\optimal{J_{\mathcal R}}\le\optimal{J_{\mathcal O}}$. Therefore
  $\optimal{J_{\mathcal R}}=\optimal{J_{\mathcal O}}$ so the solution of
  Problem~\ref{problem:rcp} is globally optimal for
  Problem~\ref{problem:ocp}~$\alev{[0,t_f]}$. \qed
\end{pf*}

Theorem~\ref{theorem:lcvx} implies that Problem~\ref{problem:ocp} is solved in
polynomial time by an SOCP solver applied to Problem~\ref{problem:rcp}. This can
be done efficiently with several numerically reliable SOCP solvers
\cite{Dueri2014}. This demonstrates that the class of \NPhard optimal control
problems defined by Problem~\ref{problem:ocp} under
Conditions~\ref{condition:observability}-\ref{condition:transversality} is in
fact of \PP complexity.

\section{Numerical Example}
\label{sec:examples}

\begin{figure}
  \centering
  \includegraphics[width=0.7\columnwidth]{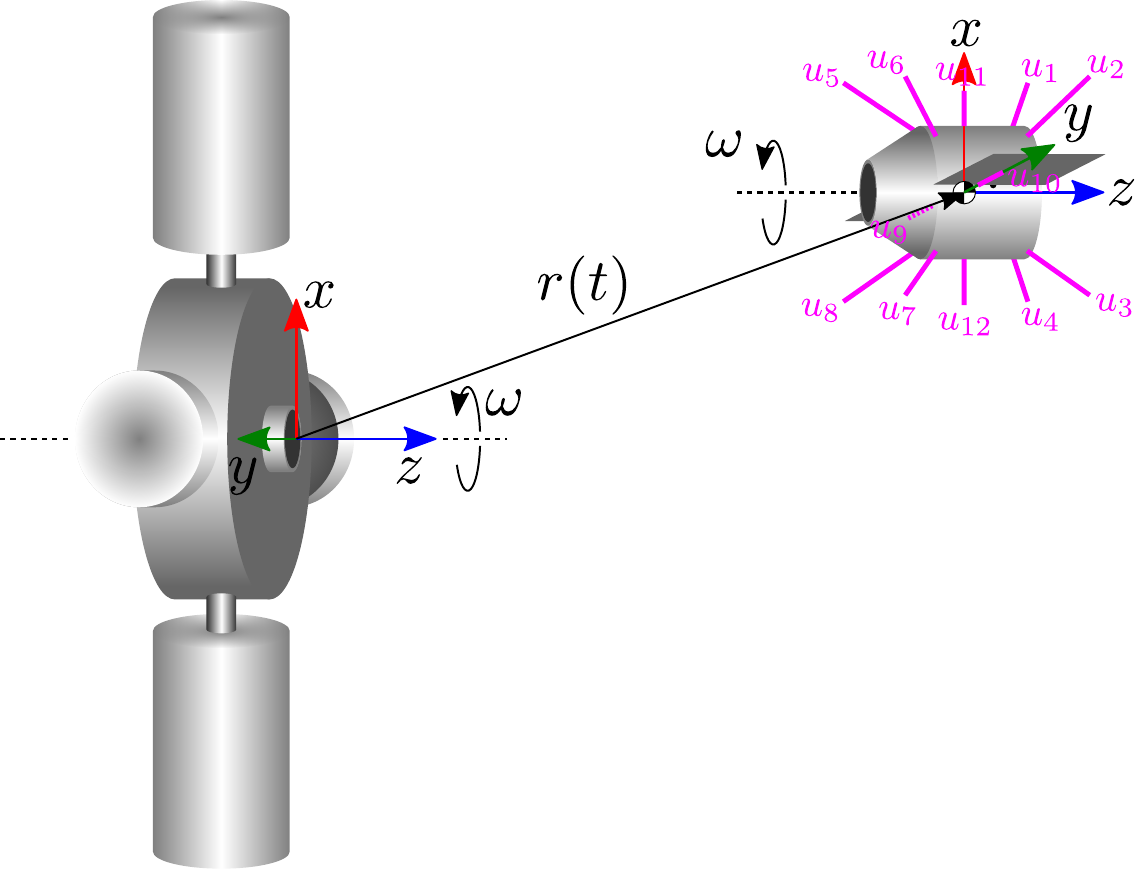}
  \caption{Spacecraft docking to a rotating space station is modelled in the
    space station's rotating frame. The spacecraft is assumed to have matched
    the space station angular velocity.}
  \label{fig:space_station_docking}
\end{figure}

This section shows how a trajectory for spacecraft docking to a rotating space
station can be computed more efficiently via Problem~\ref{problem:rcp} versus a
standard MICP approach. Python source code for this example is available
online\footnote{\url{https://github.com/dmalyuta/lcvx}}. Figure~\ref{fig:space_station_docking}
illustrates the scenario. The spacecraft's dynamics are described in the
rotating frame by:
\begin{eqnarray}
  \label{eq:dynamics}
  \dot x(t) = A(\omega)x(t)+B\inlinesum_{i=1}^Mu_i(t),
\end{eqnarray}
where $x(t)=(r(t),v(t)):\reals_+\to\reals^6$ is the position and velocity state,
$\omega\in\reals^3$ is the space station's constant angular velocity vector, and
\begin{eqnarray}
  \label{eq:AB_matrices}
  A(\omega) \definedas
  \begin{bmatrix}
    0 & I \\
    -S(\omega)^2 & -2S(\omega)
  \end{bmatrix},~%
  B\definedas
  \begin{bmatrix}
    0 \\ I
  \end{bmatrix},
\end{eqnarray}
where $S(\omega)\in\reals^{3\times 3}$ is the skew-symmetric matrix
representation of the cross product $\omega\times(\cdot)$.

\begin{figure}
  \centering
  \includegraphics[width=0.6\columnwidth]{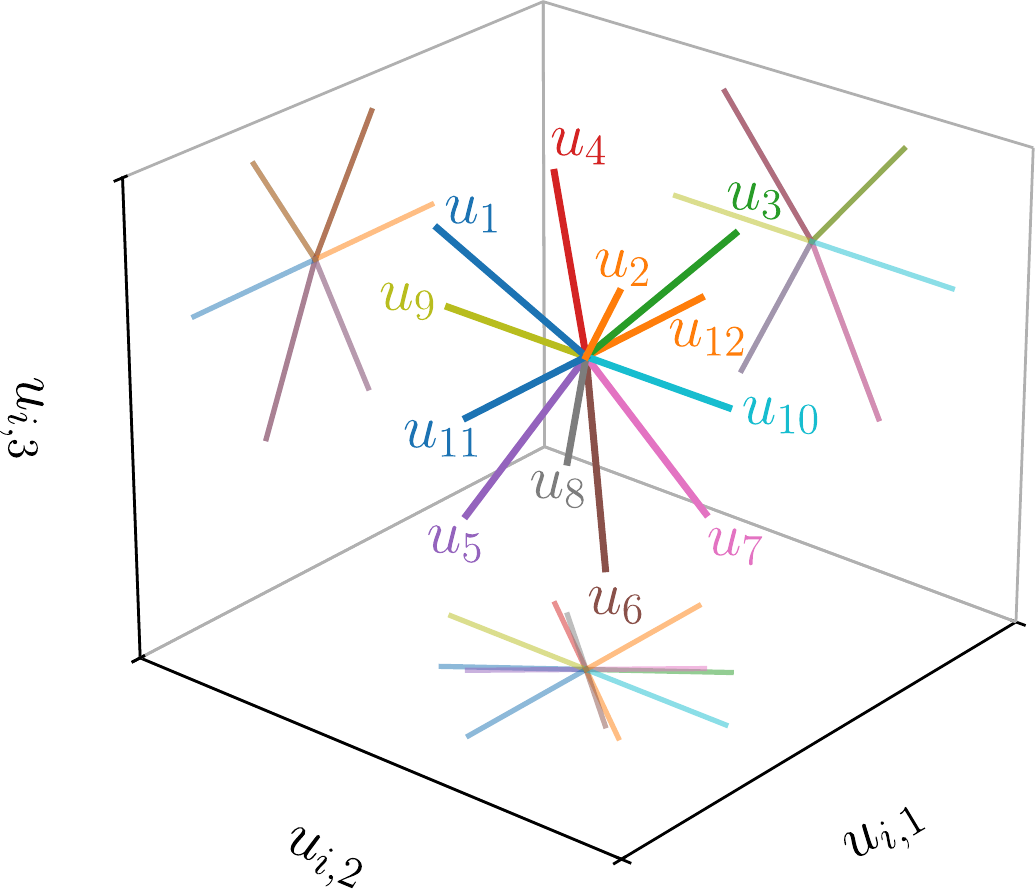}
  \caption{The spacecraft's reaction control system is capable of producing up
    to $K=4$ out of $M=12$ acceleration vectors, acting through the center of
    mass.}
  \label{fig:example_rcs}
\end{figure}

We assume that the spacecraft is equipped with a reaction control system capable
of producing up to $K=4$ acceleration vectors from a total of $M=12$ distinct
directions, as illustrated in Figure~\ref{fig:example_rcs}. Acceleration vectors
along the positive $z$-axis point with a pitch and roll of $40$ degrees. Along
the negative $z$-axis, the pitch and roll is $30$ degrees. 
The docking port is positioned at the origin and rotates with the space
station. The spacecraft is assumed to rotate with the same angular velocity
$\omega$ and is tasked to perform translation control to berth with the docking
port. We use the parameters:
\begin{gather*}
  \omega = (0,0,1)~\text{rpm},~\rho_1 = 1~\si{\milli\meter\per\second\squared},~\rho_2 = 10~\si{\milli\meter\per\second\squared}, \\
  m[t_f] = t_f,~r(0) = (5,5,100)~\si{\meter},~v(0) = (0,0,0)~\si{\meter\per\second}, \\
  r(t_f)=(0,0,0)~\si{\meter},~v(t_f)=(0,0,-0.01)~\si{\meter\per\second},
\end{gather*}
where the initial velocity choice makes the spacecraft's inertial velocity
$\omega\times r(0)$~\si{\meter\per\second}. A more realistic setting would be
$v(0)=-\omega\times r(0)$ which makes the inertial velocity zero. The motivation
for the present choice is to enrich the solution. Since the cost is linear in
$t_f$, both bisection and golden search can be applied to find the minimum $t_f$
\cite{Blackmore2010,Kochenderfer2019}.

Conditions~\ref{condition:observability}-\ref{condition:transversality} are
satisfied by this problem. Condition~\ref{condition:observability} is satisfied
since $\{-A(\omega)\T,B\T\}$ is
observable. Condition~\ref{condition:transversality} is satisfied according to
the minimum-time special case in
Figure~\ref{fig:condition_transversality_case_1}. Conditions~\ref{condition:normality}
and \ref{condition:ambiguity} hold following the discussion in
Section~\ref{subsec:condition_special_cases} and noting that the primer vector
$y(\cdot)$ in \eqref{eq:adjoint_system} can only be persistently normal to
$\omega$ or to vectors normal to $\omega$. The conditions are checked
automatically via a script in the public source code.

\begin{figure*}
  \centering
  \begin{subfigure}[b]{0.315\textwidth}
    \centering
    \includegraphics[width=0.76\textwidth]{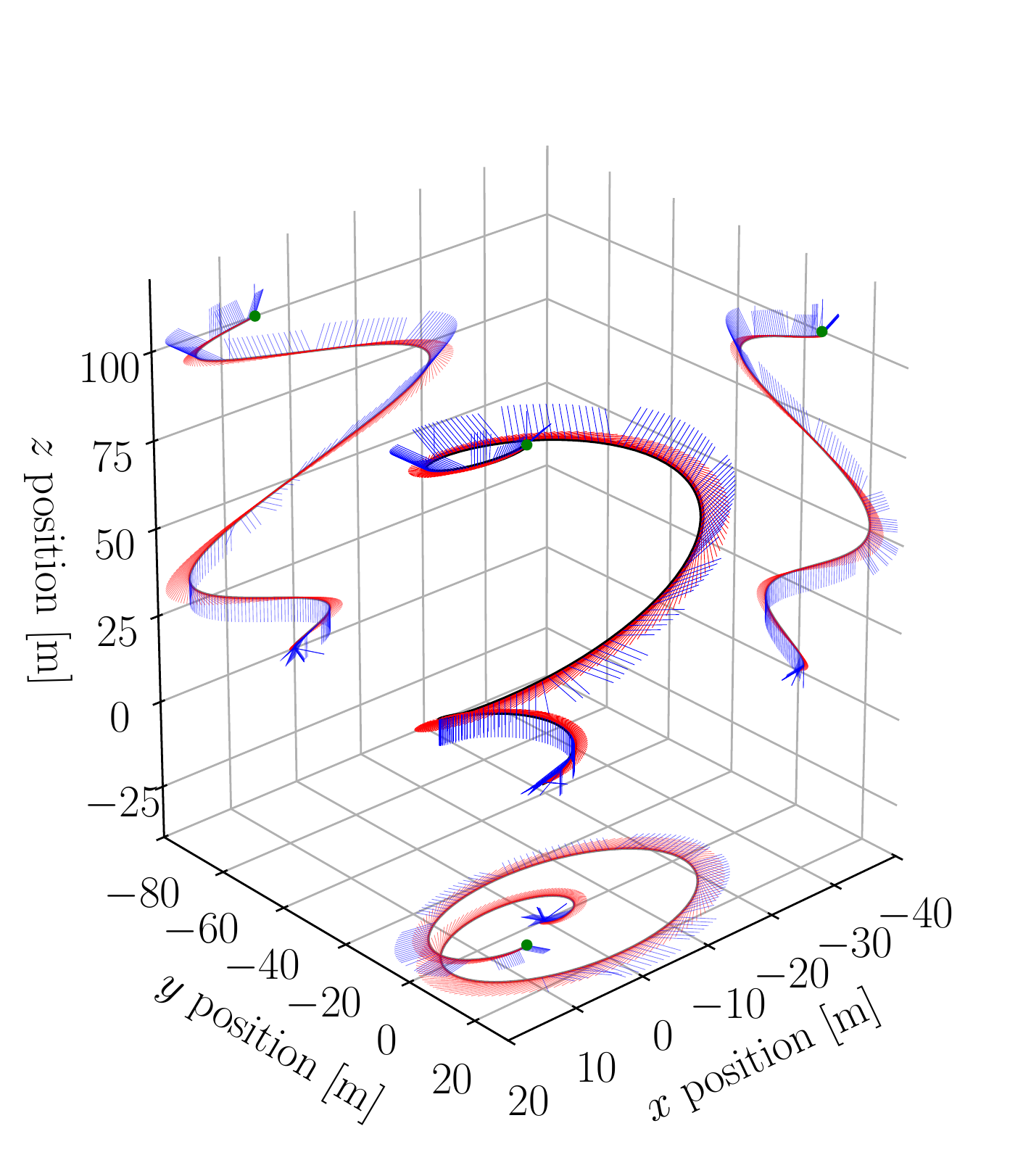}
    \vspace{-8mm}
    
    \includegraphics[width=0.76\textwidth]{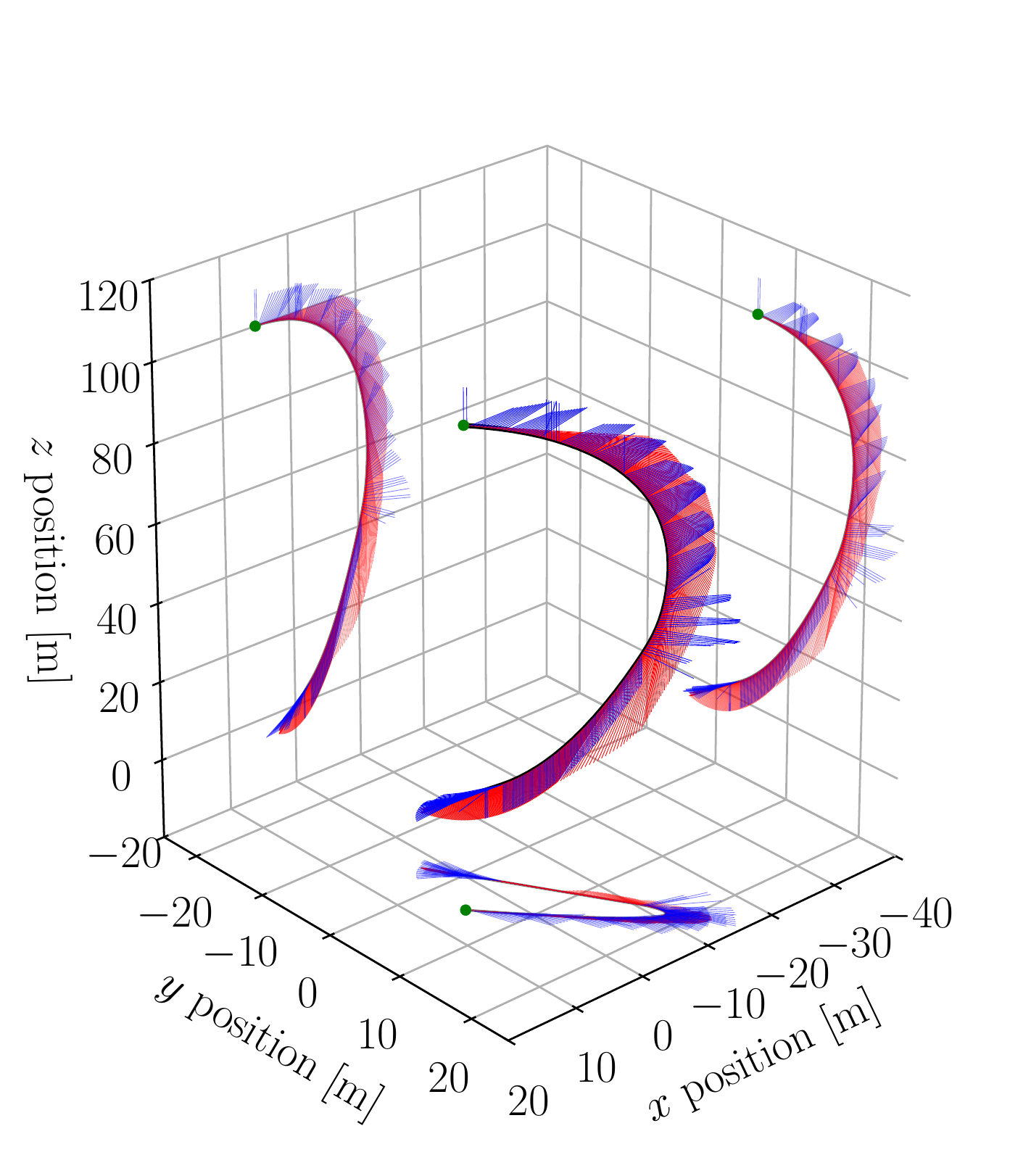}
    \caption{Trajectory in the space station rotating frame (top), and the
      inertial frame (bottom). Time of flight $t_f=135~\si{\second}$.}
    \label{fig:trajectory}
  \end{subfigure}%
  \hfill%
  \begin{subfigure}[b]{0.34\textwidth}
    \includegraphics[width=1\columnwidth]{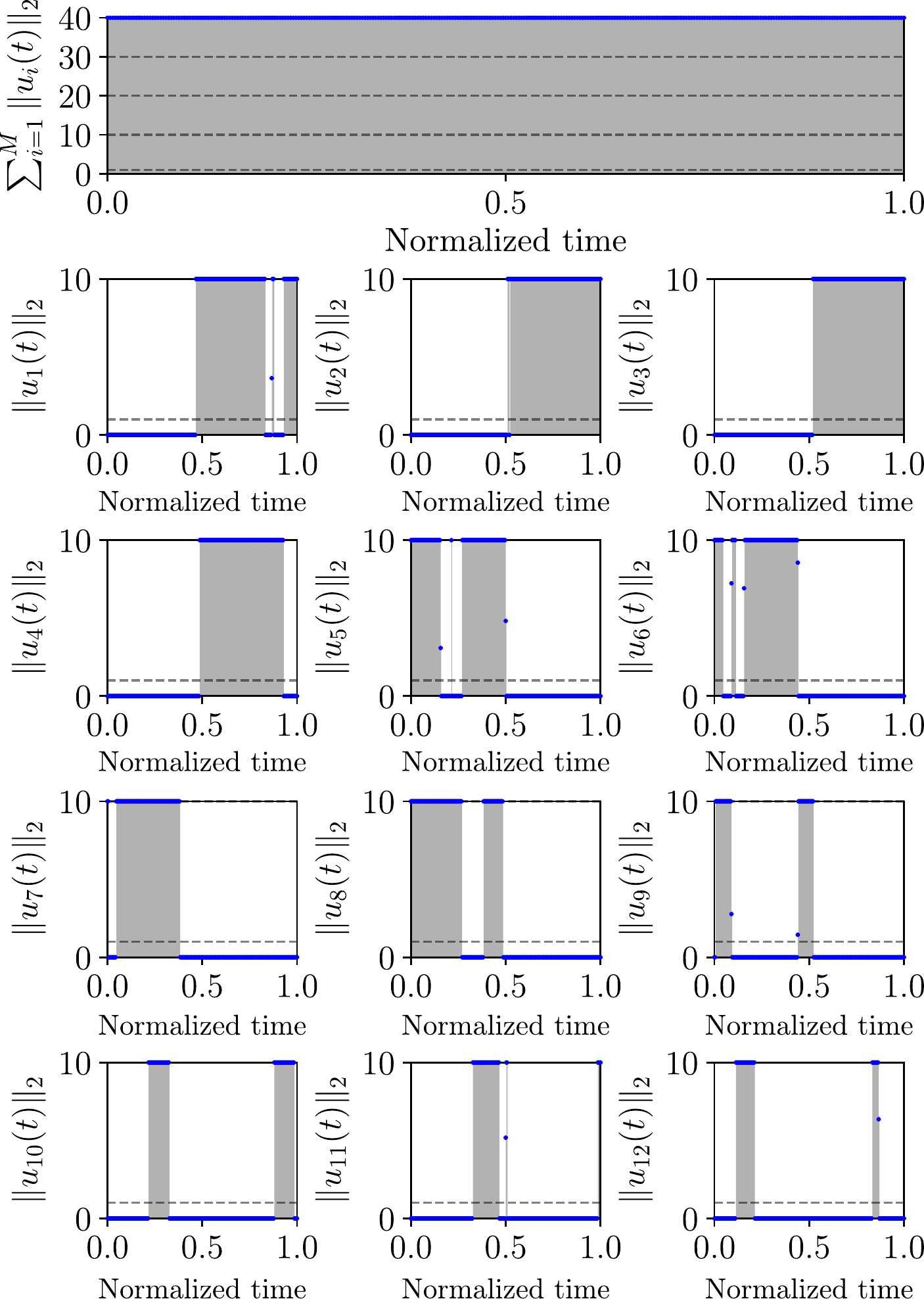}
    \caption{Optimal input norm history in
      \si{\milli\meter\per\second\squared}. The bang-bang nature of the solution
      is clearly visible. Blue markers show $\sigma_i(t)$.}
    \label{fig:input_history}
  \end{subfigure}%
  \hfill%
  \begin{subfigure}[b]{0.32\textwidth}
    \includegraphics[width=1\columnwidth]{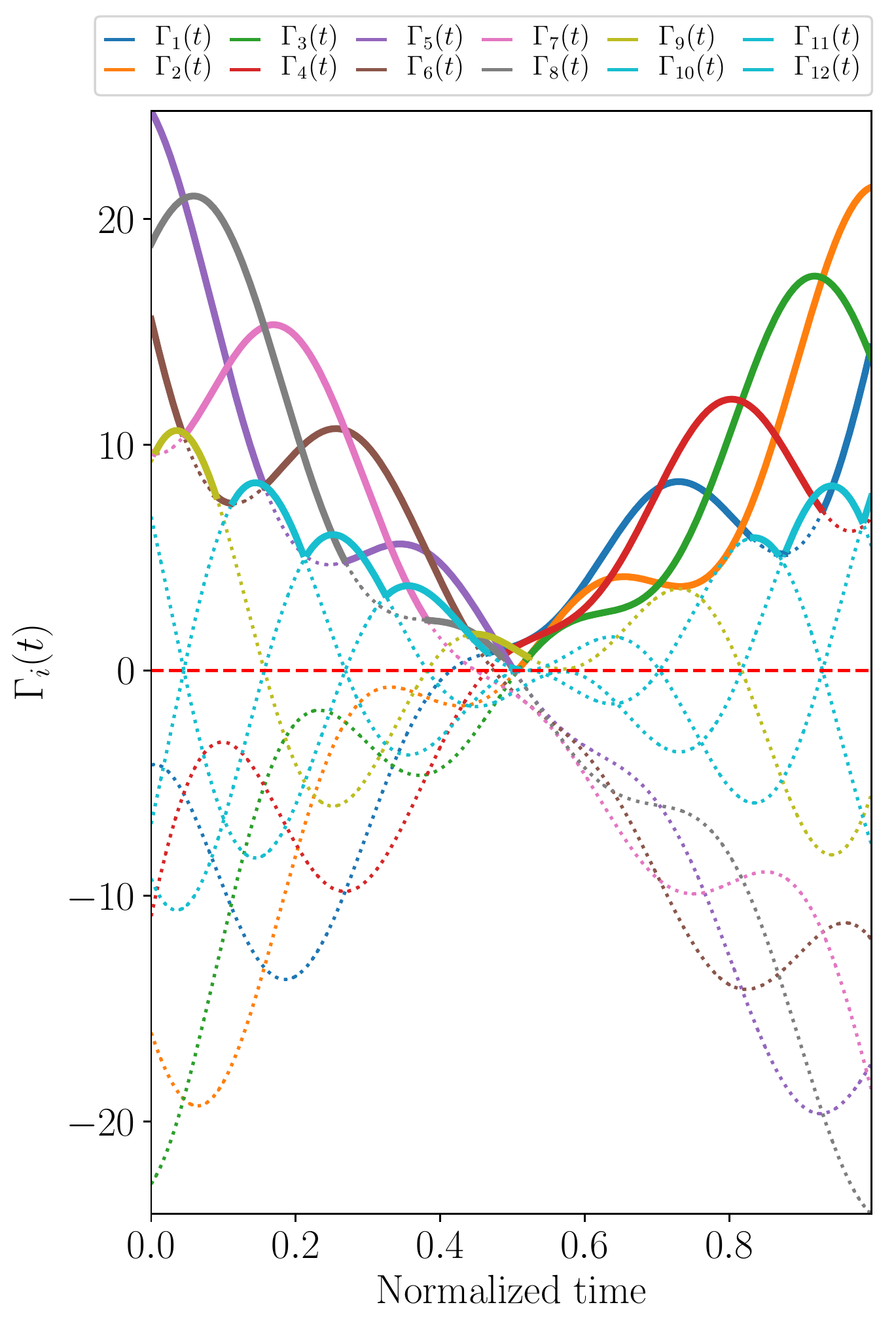}
    \caption{Time history of the input gain \eqref{eq:input_gain}. Bold lines
      show when the corresponding input is active.}
    \label{fig:input_gain}
  \end{subfigure}
  \caption{Optimal docking trajectory (\protect\subref{fig:trajectory}) together
    with the input (\protect\subref{fig:input_history}) and dual variable
    (\protect\subref{fig:input_gain}) histories. In
    (\protect\subref{fig:trajectory}), red vectors show the scaled velocity,
    blue vectors show the scaled thrust (negative of acceleration), and the
    green marker shows the initial position.}
  \label{fig:example}
\end{figure*}

The dynamics \eqref{eq:dynamics} are discretized using zeroth-order hold over a
uniform temporal grid of 300 nodes. Python~2.7.15 with ECOS~2.0.7.post1
\cite{Domahidi2013} is used for implementation on a Ubuntu~18.04.1 64-bit
platform with a 2.5~GHz Intel Core i5-7200U CPU and 8~GB of RAM. The solution
and runtime are compared to a mixed-integer formulation where \eqref{eq:ocp_d}
is implemented directly as a binary constraint using Gurobi~8.1~\cite{gurobi}.

Figure~\ref{fig:example} shows the resulting state, input and input gain
trajectories for the globally optimal solution. The solution is obtained in
$20~\si{\second}$ via Problem~\ref{problem:rcp}, whereas MICP takes an
intractable $6200~\si{\second}$. This is expected, since solving
Problem~\ref{problem:rcp} relies on an SOCP solver with polynomial time
complexity, whereas MICP has exponential time complexity.

Figure~\ref{fig:input_history} confirms that the
constraints~\eqref{eq:ocp_c}-\eqref{eq:ocp_e} are satisfied. In particular, the
thrust magnitude is bang-bang as predicted in Lemma~\ref{lemma:lcvx}. The
intermediate thrusts occuring at rising and falling edges are discretization
artifacts since the lossless convexification guarantee is only ``almost
everywhere'' in nature. These artifacts have been observed since the early days
of the lossless convexification method
\cite{Acikmese2007}. Figure~\ref{fig:input_gain} confirms the optimal input
structure \eqref{eq:gamma_optimality_structure}. In particular, for the
minimum-time solution it is always the inputs corresponding to the largest $K=4$
gain values $\Gamma_i(t)$ that are active.

\section{Future Work}
\label{sec:future_work}

The primary focus of future work is to expand the class of problems that can be
handled. This includes using different norm types and bounds in
\eqref{eq:ocp_c}, adding a lower-bound constraint to \eqref{eq:ocp_e}
(i.e. minimum number of active inputs), relaxing
Assumption~\ref{ass:pointing_set_non_overlapping} to non-overlapping input
pointing sets, considering linear time-varying dynamics in \eqref{eq:ocp_b}, and
introducing state constraints. A minor caveat of the Lemma~\ref{lemma:lcvx}
proof is that conditions which are proven to hold ``almost everywhere'' are
assumed not to fail on nowhere dense sets of positive measure (e.g. the fat
Cantor set) \cite{Morgan1990}. We do not expect this pathology to occur for any
practical problem, and in the future we seek to rigorously eliminate this
pathology.

\section{Conclusion}
\label{sec:conclusion}

This paper presented a lossless convexification method for solving a class of
optimal control problems with semi-continuous input norms. Such problems are
inherently mixed-integer in nature. By relaxing the problem to a convex one and
proving that the relaxed solution is globally optimal for the original problem,
solutions can be found via convex optimization in polynomial time. Our result
thus shows that a practical class of \NPhard optimal control problems is in fact
of \PP complexity. The algorithm is amenable to real-time onboard implementation
and can be used to accelerate design trade studies.

\section{Acknowledgements}

We thank Emma Hansen for conducting preliminary simulations.

\bibliographystyle{ieeetr}
\bibliography{references}

\end{document}